\documentclass[10pt]{amsart}
\usepackage{amsmath}
\usepackage{amsxtra}
\usepackage{amscd}
\usepackage{amsthm}
\usepackage{amsfonts}
\usepackage{amssymb}
\usepackage{eucal}

\newtheorem{theorem}{Theorem}[section]

\newtheorem{cor}[theorem]{Corollary}
\newtheorem{lem}[theorem]{Lemma}
\newtheorem{prop}[theorem]{Proposition}

\theoremstyle{definition}

\theoremstyle{remark}

\theoremstyle{remark}

\numberwithin{equation}{section}
\newcommand{\nc}{\newcommand}
\nc{\on}{\operatorname}
\nc{\ch}{\mbox{ch}}
\nc{\Q}{{\mathbb Q}}
\nc{\Z}{{\mathbb Z}}
\nc{\C}{{\mathbb C}}
\nc{\pone}{{\mathbb C}{\mathbb P}^1}
\nc{\pa}{\partial}
\nc{\F}{{\mathcal F}}
\nc{\arr}{\rightarrow}
\nc{\larr}{\longrightarrow}
\nc{\al}{\alpha}
\nc{\ri}{\rangle}
\nc{\lef}{\langle}
\nc{\W}{{\mathcal W}}
\nc{\la}{\lambda}
\nc{\ep}{\epsilon}
\nc{\su}{\widehat{{\mathfrak sl}}_2}
\nc{\sw}{{\mathfrak s}{\mathfrak l}}
\nc{\g}{{\mathfrak g}}
\nc{\h}{{\mathfrak h}}
\nc{\n}{{\mathfrak n}}
\nc{\N}{\widehat{\n}}
\nc{\G}{\widehat{\g}}
\nc{\De}{\Delta_+}
\nc{\gt}{\widetilde{\g}}
\nc{\Ga}{\Gamma}
\nc{\one}{{\mathbf 1}}
\nc{\z}{{\mathfrak Z}}
\nc{\zz}{{\mathcal Z}}
\nc{\Hh}{{\mathcal H}_\beta}
\nc{\qp}{q^{\frac{k}{2}}}
\nc{\qm}{q^{-\frac{k}{2}}}
\nc{\La}{\Lambda}
\nc{\wt}{\widetilde}
\nc{\qn}{\frac{[m]_q^2}{[2m]_q}}
\nc{\cri}{_{\on{cr}}}
\nc{\kk}{h^\vee}
\nc{\sun}{\widehat{\sw}_N}
\nc{\hh}{\widehat{\mathfrak h}}
\nc{\HH}{{\mathcal H}_{q,t}}
\nc{\ca}{\wt{{\mathcal A}}_{h,k}(\sw_2)}
\nc{\gl}{\widehat{{\mathfrak g}{\mathfrak l}}_2}
\nc{\el}{\ell}
\nc{\s}{{\mathbf s}}
\nc{\bi}{\bibitem}
\nc{\om}{\omega}
\nc{\WW}{\W_\beta}
\nc{\scr}{{\mathbf S}}
\nc{\ab}{{\mathbf a}}
\nc{\rr}{r}
\nc{\ol}{\overline}
\nc{\con}{qt^{-1} + q^{-1}t}
\nc{\den}{q^{\el-1} t^{-\el+1}+ q^{-\el+1} t^{\el-1}}
\nc{\ds}{\displaystyle}
\nc{\B}{B}
\nc{\A}{{\mathbb A}}
\nc{\GG}{{\mathcal G}}
\nc{\UU}{{\mathcal U}}
\nc{\MM}{{\mathcal M}}
\nc{\CC}{{\mathcal C}}
\nc{\GL}{{}^L G}
\nc{\dzz}{\frac{dz}{z}}
\nc{\Res}{\on{Res}}
\nc{\rep}{{\mathcal R}ep \;}
\nc{\uqg}{U_q \G}
\nc{\uqgg}{U_q \g}
\nc{\Fq}{{\mathbb F}_q}

\nc{\stimes}{\ltimes}
\nc{\K}{\hat{\mathcal K}}
\nc{\Ql}{\ol{\mathbb Q}_\ell}

\nc{\ga}{\gamma}
\nc{\PL}{{}^L P}
\nc{\E}{\mc E}
\nc{\mc}{\mathcal}
\nc{\mbf}{\mathbf}
\nc{\bb}{{\mathfrak b}}
\nc{\OO}{{\mc O}}
\nc{\Po}{{\mc P}}
\nc{\V}{{\mc V}}
\nc{\yy}{{\mc Y}}
\nc{\M}{\mathcal M}
\nc{\Coh}{{{\mathcal C}oh}}
\nc{\Cohn}{\Coh_n}
\nc{\f}{{\mathcal F}}
\nc{\si}{_E}
\nc{\Gaf}{{\mathbb G}_{a,\Fq}}
\nc{\KK}{{\mathfrak k}}
\newcommand{\Sing}{{{\rm Sing}}}
\nc{\R}{{\mathbb R}}
\newcommand{\End}{{\rm End\,}}
\newcommand{\Ve}{ L_{\bs \La}}

\newcommand{\id}{{{\rm id}}}

\newcommand{\bean}{\begin{eqnarray}}
\newcommand{\eean}{\end{eqnarray}}
\newcommand{\be}{\begin{displaymath}}
\newcommand{\ee}{\end{displaymath}}
\newcommand{\bea}{\begin{eqnarray*}}   
\newcommand{\eea}{\end{eqnarray*}}
\newcommand{\bs}{\boldsymbol}
\newcommand{\Ref}[1]{{$($\ref{#1}$)$}}

\begin{document}
\title[Counterexample to Bethe Ansatz Conjecture]
{Multiple orthogonal polynomials and \\
a counterexample to Gaudin Bethe Ansatz Conjecture}

\author{E. Mukhin and A. Varchenko}
\thanks{Research of E.M. is supported in part by NSF grant DMS-0140460.
Research of A.V. is supported in part by NSF grant DMS-0244579.}
\address{EM: Department of Mathematics, 
Indiana University-Purdue University-Indianapolis, 
402 N.Blackford St., LD 270, Indianapolis, IN 46202, USA}
\email{mukhin@math.iupui.edu}
\address{A.V.: Department of Mathematics, University of North Carolina 
at Chapel Hill, Chapel Hill, NC 27599-3250, USA }
\email{anv@email.unc.edu}
\thanks {Submitted on  05/17/2005. Revised on 10/05/2005.}.

\thanks{Mathematics Subject Classification:  82B23,\ 33C45.}

\begin{abstract}
  Jacobi polynomials are polynomials whose zeros form the unique
  solution of the Bethe Ansatz equation associated with two $sl_2$
  irreducible modules. We study sequences of $r$ polynomials whose
  zeros form the unique solution of the Bethe Ansatz equation
  associated with two highest weight $sl_{r+1}$ irreducible modules,
  with the restriction that the highest weight of one of the modules 
  is a multiple of the first fundamental weight. 
  
  We describe the recursion which can be used to compute these
  polynomials. Moreover, we show that
  the first polynomial in the sequence coincides with the 
  Jacobi-Pi\~neiro multiple orthogonal 
  polynomial and others are given by 
  Wronskian type determinants of Jacobi-Pi\~neiro polynomials.

  As a byproduct we describe a counterexample to the Bethe
  Ansatz Conjecture for the Gaudin model.
\end{abstract}
\maketitle

\section{Introduction}
\subsection{Content of the paper}
It is well known that zeros of the classical Jacobi orthogonal polynomial
with parameters $\al, \beta$ satisfy a system of algebraic equations,
which is known as the Bethe Ansatz equation of the Gaudin model
associated to $sl_2$ and two irreducible modules with highest weights
$-\al-1,-\beta-1$, see \cite{Sz}, \cite{V2}.

In this paper we 
study the Bethe Ansatz equation of the Gaudin model associated to
$sl_{r+1}$ and two irreducible modules, one of which has an arbtrary
highest weight and the highest weight of another one is a multiple of
the first fundamental weight, see \Ref{BAE}. In this case the dimensions 
of the spaces of singular vectors of a given weight are at most 
one-dimensional and therefore we may expect an appearance of some 
orthogonal polynomials. In this paper we show that this is indeed so.

We write the weights as $\bs m=(m_1,\dots,m_r)$ and $(k,0,\dots,0)$, 
see \Ref{weights}. We show that then 
for each partition 
$\bs l=(l_1,l_2\dots,l_r)$, where $l_i$ are non-negative integers, 
$l_i\geq l_j$, if $i>j$, and 
for generic values of the parameters $\bs m,k$, there exists a 
unique solution of the Bethe Ansatz equation. We consider polynomials 
$y_1,\dots,y_r$, $\deg y_i=l_i$, whose zeros form such a solution.

We call parameters $\bs m, k, \bs l$ consistent if $m_1,\dots,m_r,k$ are non-negative integers and \Ref{consistent} is satisfied.
According to \cite{MV1}, for consistent values of parameters, 
polynomials $y_1,\dots,y_r$ can be computed from a certain space of
polynomials $V(\bs m,\bs l,k)$ 
of dimension $r+1$ which has only two singular points. 
Degrees of polynomials in $V(\bs m,\bs l,k)$  and 
the exponents at the singular points 
are given explicitly in terms of the weights $\bs m, k$ 
and the partition $\bs l$, see Lemma \ref{V}.
 
We construct the spaces $V(\bs m,\bs l,k)$ recursively using explicit linear 
differential operators, see Corollary \ref{Vk op},
 and discover that $V(\bs m,\bs l,k)$ has a basis written in
terms of multiple orthogonal polynomials called the Jacobi-Pi\~neiro
polynomials, see Lemma \ref{basis lem} and Proposition \ref{v=P}.

The Jacobi-Pi\~neiro polynomial \cite{P} is the 
unique monic polynomial of degree $l_1$ 
whose coefficients are rational functions of $\bs m,k$ and 
which is orthogonal to functions
\bea
\underbrace{1,x,\dots,x^{l_1-l_2-1}}_{l_1-l_2}\ ,\underbrace{x^{-m_2-1},x^{-m_2},\dots, x^{-m_2+l_2-l_3-2}}_{l_2-l_3}\ ,
\dots, {}\hspace{40pt}\notag \\ \underbrace{x^{-\sum_{i=2}^{r}m_i-r+1},\dots,
x^{-\sum_{i=2}^{r}m_i-r+l_r}}_{l_r}
\eea
with respect to the scalar product given by
\be
(f(x),g(x))=\int_0^1f(x)g(x)(x-1)^{-k-1}x^{-m_1-1}dx,
\ee
if $m_i,k$ are all negative real numbers. 

It follows that the polynomials $y_1,\dots,y_r$ are Wronskians of the
Jacobi-Pi\~neiro polynomials up to an explicit factor. In particular
$y_1$ is exactly the Jacobi-Pi\~neiro polynomial. Using properties of 
multiple orthogonal polynomials we are able to prove that for general
values of parameters $\bs m,k$ the corresponding Bethe vector is
non-zero, see Theorem \ref{bv}.

Multiple orthogonal polynomials can be computed rather explicitly
and therefore we have explicit formulas for $y_1$. 
Similar computation allows us to obtain 
an explicit formula for $y_r$. Using these formulas, we 
find that for special values of parameters,  
the $sl_3$ Bethe Ansatz equation has no solutions and therefore 
we obtain a conterexample to the standard form of the 
Bethe Ansatz Conjecture, see Proposition \ref{counter}. 

\medskip

We also describe the spaces of polynomials $U(\bs m,\bs l,k)$ which
are dual to the spaces $V(\bs m,\bs l, k)$. These spaces $U(\bs m,\bs l,k)$ 
were
considered in \cite{Sc} (with $x$ replaced by $-x$). In \cite{Sc} the
problem of finding an explicit basis was reduced to a linear system of
equations of the size $(\sum_{i=1}^rm_i)-l_1$.

\medskip
Thus, summarizing, we establish a
relation between the Bethe Ansatz method and the theory of multiple orthogonal
polynomials. We also establish recurrence relations between
the vector spaces of polynomials of one variable, which have given
dimension $r+1$ and exactly two singular points at $x=0$ and $x=1$. We
start with a space of polynomials consisting of monomials and construct
all other spaces using suitable first order linear polynomial differential
operators. Moreover, we 
show that the differential operators form an extension of
the $sl_{r+1}$ affine Weyl group. Those differential operators provide the
Rodrigues-type explicit formulas for the bases of those spaces of polynomials.

\medskip

A study similar to the one conducted in this paper can be done
in the trigonometric case, related to multiple Laguerre polynomials,
and also in the cases associated with the $XXX$ model
and the $XXZ$ model, related to multiple little q-Jacobi polynomials studied in
\cite{PV}. See in \cite{MV3} related discussions of the $XXX$ model.

\medskip

Our paper is constructed as follows. In the rest of the
introduction we state our main results in more detail for the cases of 
$sl_2$ and $sl_3$ and present our counterexample to the standard form 
of the Bethe Ansatz Conjecture. 
In Section \ref{spaces sec} we define the spaces of polynomials
$V(\bs m,\bs l,k)$ and $U(\bs m,\bs l,k)$.  In Section \ref{recursion
  sec} we construct those spaces recursively.  In Section \ref{V sect}
we construct the basis of $V(\bs m,\bs l,k)$ in terms of
Jacobi-Pi\~neiro polynomials. In Section \ref{U sect} we give an
explicit basis of $U(\bs m,\bs l,k)$. In Section \ref{Bethe sec} we
use our results to study the Bethe Ansatz Conjecture. 
In Section \ref{app sec} we present remarks on some integral identity and
on an affine Weyl group action in our situation.

\medskip

We thank Ken McLaughlin for useful discussions. EM thanks 
Max-Planck-Institute f\"{u}r Mathematik where part of this work was done 
for hospitality and a stimulating enviroment.

\subsection{%Main results in the case of $sl_2$}
Jacobi polynomials with negative integer parameters}
\label{subsection r=1}
Let $m,l,k$ be non-negative integers such that
\bean\label{1st ineq}
k \geq l \geq 0 ,
\qquad
m \geq l \geq 0 .
\eean
Let
\bean
\mc D(l,k)\ =\ 
x(x-1)
\frac{d^2}{d x^2}\ -\
(k x + m(x-1))
\frac{d}{d x}\ +\
l(k+m+1-l) \ .
\notag
\eean
Denote by $V(l,k)$ the two-dimensional complex vector
space of solutions of the differential
equation $\mc D(l,k) f = 0$. The space $V(l,k)$
has the properties:
\begin{enumerate}
\item[(i)]
The space $V(l,k)$ consists of polynomials
of degrees $l$ and $k+m+1-l$, \linebreak
$l < k+m+1-l$.

\item[(ii)] There exists a polynomial in $V(l,k)$
which has a root at 0 of multiplicity $m+1$.

\item[(iii)] There exists a polynomial in $V(l,k)$
which has a root at 1 of multiplicity $k+1$.
\end{enumerate}
These three properties define the two-dimensional space $V(l,k)$ uniquely.

\bigskip

{\it Example.} Polynomials 1 and $x^{m+1}$ form a basis in $V(0,0)$.
\hfill
$\square$

\bigskip

Let $y(l,k;x) $ be the monic polynomial in $V(l,k)$ of degree $l$.
Then 
$$
y(l,k;x)\ =\ P_l^{(\al,\beta)}(x)\ ,
$$ 
where $\al=-k-1$, $\beta=-m-1$, 
and $P_l^{(\al,\beta)}(x)$ is the  monic Jacobi polynomial on the interval
$[0,1]$.

A basis in  $V(l,k)$ is formed by $P_l^{(-k-1,-m-1)}(x)$ and
$x^{m+1} P_{k-l}^{(-k-1,m+1)}(x)$.

According to the Rodrigues formula the polynomial $y(l,k;x)$ is equal to
\bea
x^{m+1} (x-1)^{k+1} \ \frac{d^l}{dx^l} \
[\ x^{l-m-1} (x-1)^{l-k-1}\ ]
\eea
up to multiplication by a non-zero constant.

\bigskip

\noindent
Consider the first order linear differential operators
\begin{align}
&
D_0(l,k)\ =\ x(x-1)\frac d {dx}\ -\ l (x-1)\ - \ k - 1\ , 
\notag
\\
&
D_1(l,k)\ =\ x(x-1)\frac d {dx}\ -\ (k + m + 1 - l)(x-1)\ -\ k-1\ ,
\notag
\\
&
D^\vee_0(l,k)\ =\ x\frac d {dx}\ -\ (k + m + 1 - l) \ ,
\notag
\\
&
D^\vee_1(l,k)\ =\ 
  x \frac d {dx}\ -\ l\ . 
\notag
\end{align}
The operators commute:
\bea
D_1(l,k+1)\ D_0(l,k)\ &=&\ D_0(l+1,k+1)\ D_1(l,k)\ ,
\\
D_1^\vee(l,k-1)\ D_0^\vee(l,k)\ &=&\
D_0^\vee(l-1,k-1) \ D_1^\vee(l,k)\ .
\eea
The operators define isomorphisms
\bea\label{maps r=1}
%& &
%\\
D_0(l,k)& :& V(l,k)\ \to\ V(l,k+1)\ ,
\notag
\\
D_1(l,k)& :& V(l,k) \ \to \ V(l+1,k+1)\ ,
\notag
\\
D_0^\vee(l,k)& :& V(l,k)\ \to\ V(l,k-1)\ ,
\notag
\\
D_1^\vee(l,k)& :& V(l,k)\ \to \ V(l-1,k-1)\ ,
\notag
\eea
if the preimage and image spaces are defined, that is if
the corresponding parameters satisfy \Ref{1st ineq}. 
In particular, the composition
\bea
D_1(l-1,k-1) \cdots D_1(0,k-l) \ D_0(0,k-l-1) \cdots D_0(0,0)
\eea
defines an isomorphism $V(0,0) \to V(l,k)$.

%Each of the isomorphisms %\Ref{maps r=1}
%sends the triple of lines
%$L_\infty, L_0, L_1$ of the preimage space to the corresponding triple of the
%image space. 
We have
\bean\label{recur r=1}
%&
%\\
y(l,k+1;x) \ &=& \ \frac{k+m+2-2l}{(k+m+2-l)(k+1-l)}\
D_0(l,k) \ y(l,k;x) \ ,
%\notag
\\
y(l+1,k+1;x) \ &=& \ -\ \frac1{k+m+1-2l}\
D_1(l,k) \ y(l,k;x) \ .
\notag
\eean

%The compositions 
%\bea
%D^\vee_0(l,k+1)\ D_0(l,k)\ , &&
%D_0(l,k-1) \ D^\vee_0(l,k)\ , 
%\\
%D^\vee_1(l+1,k+1)\ D_1(l,k)\ ,&&
%D_1(l-1,k-1) \ D^\vee_1(l,k) 
%\eea
% define scalar operators on $V(l,k)$ equal to 
%\bean\label{constatants}
%&&
%\\
%(l-k-1) (m+k+2-l)\ ,
%&&
%(l-k) (m+k+1-l) \ , 
%\notag
%\\
%(l+1) (m-l)\ , &&
%l (m+1-l)  \ ,\
%\notag
%\eean
%respectively.

\bigskip

For any $i$, denote by $ W(f_1, \dots ,f_i)$ the Wronskian of functions 
$f_1, \dots , f_i$ of $x$.

The formula 
\bea
\langle f,g \rangle_{l,k} \ = \ \frac{W(f,g)}{x^k (x - 1)^m}\ ,
\qquad
f,g \in V(l,k)\ ,
\eea
defines on $V(l,k)$ a non-degenerate skew-symmetric bilinear form.

The operators $D_0(l,k)$ and $\ D_1(l,k)$ are adjoint to the operators
$-D^\vee_0(l,k+1)$ and $-D^\vee_1(l+1,k+1)$  respectively in the following
sense:
\begin{align*}\label{Adj}
\langle D_0(l,k) f  ,\ g \rangle_{l,k+1} &= - 
\langle f  ,\ D^\vee_0(l,k+1) g \rangle_{l,k},
\qquad f \in V(l,k), g \in V(l,k+1) ,
\notag
\\
\langle D_1(l,k) f , g \rangle_{l+1,k+1}& = - 
\langle f , D^\vee_1(l+1,k+1) g \rangle_{l,k},\ 
f \in V(l,k), g \in V(l+1,k+1) .
\notag
\end{align*}

\bigskip

Let $t_1, \dots , t_l$ be the roots of the polynomial $y(l,k;x)$. The roots are 
pairwise distinct, different from 0 and 1, and form a solution to the
system of equations
\bean\label{BAE r=1}
\frac{m}{t_i}\ +\ \frac{k}{t_i-1}\ -\
\sum_{j,\ j\neq i} \frac2{t_i - t_j}\ =\ 0 \ ,
\qquad i = 1, \dots , l .
\eean
If parameters $m, l, k$ satisfy \Ref{1st ineq}, then system \Ref{BAE r=1}
has a unique solution $t_1, \dots , t_l$ up to permutations of these numbers.

%System \Ref{BAE r=1} depends on parameters $m, l, k$ satisfying \Ref{1st ineq}. 
Formulas \Ref{recur r=1} give a recurrent way to generate solutions to 
\Ref{BAE r=1}.

\bigskip 

Coefficients of the polynomial $y(l,k;x)$ 
are rational functions of $k,  m$.
%which are regular if Re $k < -1$ and Re $m < -1$. The polynomial $y(l,k;x)$
%is well defined for those values of parameters.

Assume that Re $k < 0$ and Re $m < -1$.
For polynomials $f, g$, \
set
\bea
(f,\ g)_{k}\ = \ \int_0^1 \
 \frac{f(x)\ g(x)}{ x^{m+2}\ (1-x)^{k+1}} \ dx\ .
\eea
The operators $D_0(l,k)$ and $\ D_1(l,k)$ are adjoint to the operators
$-D^\vee_0(l,k+1)$ and $-D^\vee_1(l+1,k+1)$ respectively in the following
sense:
\bea
( D_0(l,k) f  ,\ g )_{k+1} &= &- 
( f  ,\ D^\vee_0(l,k+1) g )_{k}, \\
( D_1(l,k) f ,\ g )_{k+1} &= &- 
( f ,\ D^\vee_1(l+1,k+1) g )_{k}  .
\notag
\eea
We have
\bea\label{orth r=1}
(y(l,k;x), \ x^n)_{k}\ = \ 0\ ,
\qquad
1 \leq n \leq l \ ,
\eea
and
\bea\label{norms r=1}
(y(l,k;x), \ y(l,k;x))_{k}\ = \ l!\
\prod_{i=1}^l
\frac {m+1-i}
{(k+m+2-l-i)^2}\ (1, 1)_{k-l}\ ,
\eea
where
\bea
(1,1)_{k-l}\ =\
\frac
{  \Gamma(-m-1)\ \Gamma(-k+l)}
{  \Gamma(-k-m-1+l) }\ .
\eea

%\bigskip

%\noindent
%{\bf Remark on two $sl_2$ actions.} 
% The differential operators satisfy the following relations:
% \bea
%D_1^\vee(l,k+1)\ D_0(l,k)\ &=&\
%D_0(l-1,k-1) \ D_1^\vee(l,k)\ .
%\\
%D_0^\vee(l+1,k+1)\ D_1(l,k)\ &=&\
%D_1(l,k-1) \ D_0^\vee(l,k)\ ,
%\\
%D^\vee_0(l,k+1)\ D_0(l,k)\  + \ m\ +\ 2 + 2 (k-l)\  &=& \  
%D_0(l,k-1) \ D^\vee_0(l,k)\ ,
%\\
%D^\vee_1(l+1,k+1)\ D_1(l,k)\ +\ 2l - m\  & = &  \
%D_1(l-1,k-1) \ D^\vee_1(l,k)\  .
%\eea
% Let
%\bea
%W \ =\ \oplus_{l, k \in \Z} W(l,k)\ ,
%\qquad W(l,k) = \C[x] \ .
%\eea
%Define operators $E_i, F_i, H_i : W \to W$, $i = 0, 1$. For $w \in W_{l,k}$,
%set 
%\bea
%&&
%E_0 w = D^\vee_0(l,k) w ,
%\qquad
%F_0 w = D_0(l,k) w ,
%\qquad
%H_0 w = (-m-2-2k+2l) w ,
%\\
%&&
%E_1 w = D^\vee_1(l,k) w ,
%\qquad
%F_1 w = D_1(l,k) w ,
%\qquad
%H_1 w = (m-2k)  w . 
%\phantom{aaaaa}
%\eea
% Then every 
%generator with the subscript 0 commutes with every generator with
%the subscript 1, and
%\bea
%[E_i,F_i]\ =\ H_i\ , 
%\qquad
%[H_i,E_i]\ =\ 2 E_i\ , 
%\qquad
%[H_i,F_i]\ =\ - \ 2 F_i\ . 
%\eea

\subsection{%Main results in the case of $sl_3$}
Three-dimensional spaces of polynomials with two singular points}
\label{subsection r=2}
The parameters $ \bs m = (m_1, m_2), 
\ \bs l = (l_1, l_2) \in \Z_{\geq 0}^2, \ k \in \Z_{\geq 0}$
are called {\it consistent} if
\bean\label{2d ineq}
k\ \geq l_1\ \geq \ l_2\ \geq\ 0 \ ,
\qquad
m_1 \ \geq\ l_1 - l_2 \ , 
\qquad
m_2\ \geq\ l_2\ . 
\eean
Let $ \bs m = (m_1, m_2), \ \bs l = (l_1, l_2), \ k$ be consistent. Let
\bea
\mc D(\bs l, k)\  =\ x(x-1)\frac{d^3}{d x^3}\ 
-\ ((2m_1+m_2)(x-1) + 2kx)
\frac{d^2}{d x^2}\ +\ v_1\frac{d}{d x}
\ +\ v_0\ ,
\eea
where
\bea
&&v_1 = - l_2^2+l_2(m_2+1)+l_1(k+m_1+1-l_1+l_2)+\\ && \hspace{2.2cm} (m_1+k)(m_1+m_2+k+1) 
-\frac{m_1(m_1+m_2+1)}{x}+\frac{k(k+1)}{x-1},\\ 
&&  v_0 = - \frac{(k+1)(-l_2^2+l_2(m_2+1)+l_1(m_1+k+1-l_1+l_2))}{x(x-1)}- \\
&&\hspace{4 cm}\frac{l_1(l_2-l_1+m_1+k+1)(2-l_2+m_1+k+m_2)}{x} .\\
\eea
Denote by $V(\bs l,k)$ the three-dimensional complex vector
space of solutions of the differential
equation $\mc D(\bs l, k) f = 0$. The space $V(\bs l,k)$
has the properties:
\begin{enumerate}
\item[(i)]
The space $V(\bs l,k)$ consists of polynomials
of degrees $l_1$, $m_1+1 + k + l_2 - l_1$,  and $m_1+m_2 + 2 + k - l_2$, 
where
$$
l_1\ < \ m_1+1+k+l_2-l_1\ <\ m_1+m_2+2+k-l_2\ .
$$
\item[(ii)]   There exist three polynomials in
$ V(\bs l,k)$ which have a root at 0 of multiplicities
$0,\ m_1+1,\ m_1+m_2+2$, respectively.
\item[(iii)]
There exist three polynomials in
$ V(\bs l,k)$ which have a root at 1 of multiplicities
$0,\ k+1,\ k+2$,\ respectively.

\end{enumerate}
These three
properties define the three-dimensional space $V(\bs l,k)$ uniquely.

\bigskip

{\it Example.} Polynomials 1, $x^{m_1+1}$,  and $x^{m_1+m_2 + 2}$
form a basis in $V(\bs 0,0)$, where $\bs 0 = (0,0)$.
\hfill
$\square$

\bigskip

Set
\bea
e_0(\bs l,k) = l_1 ,
\quad
 e_1(\bs l,k)  =   m_1+1+k+l_2-l_1 ,
\quad
 e_2(\bs l,k)  =  m_1+m_2+2+k-l_2 .
\eea
Set  
$$
T_1(x) \ =\  x^{m_1} (x-1)^{k}\ ,
\qquad
 T_2(x)\ = \ x^{m_2}\ .
$$ 
Let $f_{0}, f_{1}$ be the monic polynomials in 
$V(\bs l,k)$ of degrees $e_0(\bs l,k), \ e_1(\bs l,k)$,  respectively.
Set
\bean\label{y_1,y_2}
y_1(\bs l, k; x) = f_{0}(x)\ ,
\qquad
y_2(\bs l, k; x) = 
\frac {W( f_{0} , f_{ 1})}
{(e_0(\bs l,k)-e_1(\bs l,k))\ T_1(x)} \ .
\eean
These are monic polynomials of degrees $l_1$ and $l_2$, respectively.

\bigskip

 Let
\bea
\mc D^\vee(\bs l, k)\  =\ x(x - 1) \frac{d^3}{d x^3}\ -\
((2m_2 + m_1)(x-1) + k x) \frac{d^2}{d x^2}\
+ \ u_1 \frac{d}{d x}\ +\ u_0\ ,
\eea
where
\bea
&&u_1\ =\ l_2(l_1-l_2+m_2+1)+(l_2+m_2+1)(m_2+m_1+k+2-l_2)-\\
&&\hspace{4.5cm}-2-2m_2-m_1-k 
\ -\ \frac{m_2(m_1+m_2+1)}{x}\ ,
\\
&&
u_0\ =\ \frac{l_2(l_1-l_2+m_2+1)(-2-m_2-m_1-k+l_1)}{x}\ .
\eea
Denote by $U(\bs l,k)$ the three-dimensional complex vector
space of solutions of the differential
equation $\mc D^\vee(\bs l, k) f = 0$. The space $U(\bs l,k)$
has the properties:
\begin{enumerate}
\item[(i)]
The space $U(\bs l,k)$ consists of polynomials
of degrees $l_2$, $m_2 + 1 + l_1 - l_2$,  and $m_1 + m_2 + 2 + k - l_1$, 
where
$$
l_2\ < \ m_2 + 1 + l_1 - l_2 \  < \  m_1 + m_2 + 2 + k - l_1\ .
$$
\item[(ii)]   There exist three polynomials in
$ U(\bs l,k)$ which have a root at 0 of multiplicities
$0,\ m_2 + 1, \ m_1 + m_2 + 2$, respectively.

\item[(iii)]
 There exist three polynomials in
$ U(\bs l,k)$ which have a root at 1 of multiplicities
$0, \ 1,\ k+2$,\ respectively.

\end{enumerate}
These three
properties define the three-dimensional space $U(\bs l,k)$ uniquely.

\bigskip

{\it Example.} Polynomials 1, $x^{m_2+1}$,  and $x^{m_1+m_2 + 2}$
form a basis in $U(\bs 0,0)$.
\hfill
$\square$

\bigskip

\noindent
Set
\bea
e_0^\vee(\bs l,k) = l_2,
\qquad
e_1^\vee(\bs l,k) =  m_2 + 1 + l_1 - l_2,
\qquad
e_2^\vee(\bs l,k) =  m_1 + m_2 + 2 + k - l_1.
\eea
For $f_1, f_2, f_3 \in V(\bs l,k)$, set
\bea
W^\dagger_V (f_1,f_2)\ = \ \frac{W(f_1,f_2)}{T_1}\ ,
\qquad
W^\dagger_V (f_1,f_2,f_3)\ = \ \frac{W(f_1,f_2)}{T_1^2\ T_2}\ .
\eea
For given $f_1,f_2,f_3$, the function $W^\dagger_V (f_1,f_2)$
is a polynomial, and
the function $W^\dagger_V (f_1,f_2,f_3)$ is constant. 
The map $(f_1,f_2,f_3) \mapsto W^\dagger_V (f_1,f_2,f_3)$
defines a volume form on $V(\bs l,k)$.

For $g_1, g_2, g_3 \in U(\bs l,k)$, set
\bea
W^\dagger_U (g_1,g_2)\ = \ \frac{W(g_1,g_2)}{T_2}\ ,
\qquad
W^\dagger_U (g_1,g_2,g_3)\ = \ \frac{W(g_1,g_2)}{T_2^2 T_1}\ .
\eea
For given $g_1,g_2,g_3$, the function $W^\dagger_U (g_1,g_2)$
is a polynomial, and the function $W^\dagger_U(g_1,g_2,g_3)$ is constant. 
The map $(g_1,g_2,g_3) \mapsto W^\dagger_U (g_1,g_2,g_3)$
defines a volume form on $U(\bs l,k)$.

We have
\bea
U(\bs l,k)&=&  \{ W^\dagger_V(f_1,f_2)\ |\ f_1,f_2 \in V(\bs l,k) \}\ ,
\\
V(\bs l,k)&=& \{ W^\dagger_U(g_1,g_2)\ |\ g_1,g_2 \in U(\bs l,k) \}\ .
\eea
In particular, the polynomial $y_2(\bs l,k;x)$ defined in \Ref{y_1,y_2} is the unique monic
polynomial in $U(\bs l,k)$ of degree $l_2$. 

The spaces $V(\bs l,k)$ and $U(\bs l,k)$ are dual with respect 
to the following pairing:
\bea
\langle f,\ g\rangle_{\bs l,k} =\ W_V^\dagger(f, f_1,f_2)\ , 
\qquad {\rm if} \ \ g = W_V^\dagger(f_1,f_2)\ .
\eea
In particular, we have 
$$
\langle\ y_1(\bs l,k;x),\ y_2(\bs l,k;x)\ \rangle_{\bs l,k}\ =\ 0\ .
$$

\bigskip 

\noindent
Consider the first order linear differential operators
\bea
D_i(\bs l,k)& = &x(x-1)\frac d {dx}\ -\ e_i(\bs l,k) (x-1)\ -\ k - 1\ , 
\\
D_i^\vee(\bs l,k)& =& x\frac d {dx}\ -\ e_{2-i}^\vee(\bs l,k) \ ,
\eea
$i = 0, 1, 2$.
There are linear relations:
\bean\label{lin rel}
A_{12}(\bs l,k) \ D_0(\bs l,k)\ +\ A_{20}(\bs l,k)\ D_1(\bs l,k)\ +\
A_{01}(\bs l,k) \ D_2(\bs l,k)\ =\ 0\ ,
\notag
\\
A_{12}(\bs l,k) \ D_0^\vee(\bs l,k)\ +
\ A_{20}(\bs l,k)\ D_1^\vee(\bs l,k)\ +\
A_{01}(\bs l,k) \ D_2^\vee(\bs l,k)\ =\ 0\ ,
\notag
\eean
where $A_{ij}(\bs l,k) = e_i(\bs l,k) - e_j(\bs l,k)$.
Set
\bea
\bs 1_0 \ =\ (0,0)\ ,
\qquad
\bs 1_1 \ =\ (1,0)\ ,
\qquad
\bs 1_2 \ =\ (1,1)\ .
\eea
The operators commute:
\bea
D_i(\bs l+\bs 1_j,k+1)\ D_j(\bs l,k)\ &=&\ D_j(\bs l+\bs 1_i,k+1)\ D_i(\bs l,k)\ ,
\\
D_i^\vee(\bs l-\bs 1_j,k-1)\ D_j^\vee(\bs l,k)\ &=&\
D_j^\vee(\bs l-\bs 1_i,k-1) \ D_i^\vee(\bs l,k)\ 
\eea
for all $i, j$.

The operators define isomorphisms
\bea\label{maps r=2}
%& &
%\\
D_i(\bs l,k)& :& V(\bs l,k)\ \to\ V(\bs l+\bs 1_i,k+1)\ ,
\notag
\\
D^\vee_i(\bs l,k)& :& U(\bs l,k) \ \to \ U(\bs l-\bs 1_i,k-1)\ ,
\notag
\eea
for $i = 0, 1, 2$, if the preimage and image spaces are defined, that is if
the corresponding parameters satisfy \Ref{2d ineq}. 
In particular, the composition
\bea
 D_2((l_1-1,l_2-1),k-1)\  \cdots\  D_2((l_1-l_2,0),k-l_2) \ D_1((l_1-l_2-1,0),k-l_2-1)
\\
\hspace{2cm}\cdots\ D_1((0,0),k-l_1) \  D_0((0,0),k-l_1-1)\ \cdots\ D_0((0,0),0)
\eea
defines an isomorphism $V(\bs 0,0) \to V(\bs l,k)$,
and the composition
\bea
D_0^\vee((0,0),1)\  \cdots\  
D_0^\vee((0,0),k-l_1) \ D_1^\vee((1,0),k-l_1+1)
\phantom{aaaaaaaaaaaaaaaaaaa}
&&
\\
\phantom{aaaaaa}
 \cdots\ D_1^\vee((l_1-l_2,0),k-l_2) \  D_2^\vee((l_1-l_2+1,1),k-l_2+1)\ 
\cdots\ D_2^\vee((l_1,l_2),k) \phantom{}
&&
\eea
defines an isomorphism $U(\bs l,k) \to U(\bs 0,0)$.

Set $l_0 = k$,\ $l_3=0$, \ and
\bea
A_0(\bs l,k)\ =\ {}
\prod_{i=0}^2 \ {}
\frac
{k + \sum_{j=1}^i m_j - l_1 + i + 1}
{k + \sum_{j=1}^i m_j - l_1 + i + 1 - l_i + l_{i+1}}\ ,
\eea
\bea
&&A_1(\bs l,k)\ =\ 2 l_1 - l_2 - k - m_1 - 1 ,\hspace{0.7cm}
A_2(\bs l,k)\ =\ l_1 + l_2 - k - m_1 - m_2 - 2,
\\
&& A_0^\vee(\bs l,k)\ =\ l_1 + l_2 - k - m_1 - m_2 - 2, \ 
A_1^\vee( \bs l,k)\ =\ 2l_2 - l_1 - m_2 - 1,
\eea
\bea
A^\vee_2(\bs l,k)\ = \ -
\prod_{i=0}^2 \ {}
\frac
{\sum_{j=3-i}^2 m_j - l_2 + i}
{\sum_{j=3-i}^2 m_j - l_2 + i + 1 - l_{2-i} + l_{3-i}}\ .
\eea
Then for $i = 0, 1, 2,$ we have
\bea\label{2 recc}
A_i(\bs l,k)\ y_1(\bs l+\bs 1_i,k+1;x)\ &=&\ D_i(\bs l,k)\
y_1(\bs l,k;x)\ ,
\notag
\\
A_i^\vee(\bs l,k)
\ y_2(\bs l-\bs 1_i,k-1;x) &=&\  D_i^\vee (\bs l,k)\ y_2(\bs l,k;x)
\ .
\notag
\eea
%In particular, up to multiplication by a nonzero constant,
%$y_1(\bs l,k;x)$ is the image of $1$ under the operator
%\bea
%D_2((l_1-1,l_2-1),k-1)\  \cdots\  D_2((l_1-l_2,0),k-l_2) \ \times
%\phantom{aaaaa}
%&&
%\\
%\phantom{aaaaaaaa}
% D_1((l_1-l_2-1,0),k-l_2-1) \cdots\ D_1((0,0),k-l_1) \ .
%\eea
Here is a Rodrigues type formula for $y_1(\bs l,k;x)$: 
\bea
x^{m_1+m_2+2}(x-1)^{k+1}
\left( \frac{d}{dx}\right)^{l_2} \left[
x^{-m_2+l_2-1}
\left(\frac{d}{dx}\right)^{l_1-l_2} \left[
x^{l_1-l_2-m_1-1}(x-1)^{l_1-k-1}\right] \right]
&&
\\
=\ c\ y_1(\bs l,k;x)\ ,
&&
\eea
where $c$ is a non-zero constant.
%\bea
%C\ =\ (-1)^{l_1}\ {}
%\prod_{i=1}^{l_1-l_2} (k+m_1+2-l_1-i)\ {}
%\prod_{i=1}^{l_2} (k+m_1+m_2+3-l_1-i)\ {} .
%\eea
If
\bea
y_1(\bs l,k;x) \ = \ \sum_{n=0}^{l_1}\ a_{n}(\bs l,k)\, (x-1)^n \ ,
\qquad
y_2(\bs l,k;x) \ = \ \sum_{n=0}^{l_2}\ b_{n}(\bs l,k)\, x^n \ 
\eea
with $a_{l_1}(\bs l,k) = 1,\ a_{-1}(\bs l,k) = a_{l_1+1}(\bs l,k) = 0$, 
and $b_{l_2}(\bs l,k) = 1$,\ $b_{l_2+1}(\bs l,k) = 0,$ \ then
\bean\label{2' recc}
&&
 a_{n}(\bs l+\bs 1_i,k+1)\ =\ 
-\frac {k + 1 + n}{A_i(\bs l,k)}\ a_{n}(\bs l,k)\
+\ 
\frac{n - 1 - e_i(\bs l,k)}{A_i(\bs l,k)}\ a_{n-1}(\bs l,k) \ ,
\notag
\\
&&
 b_n (\bs l+\bs 1_i,k+1)\ =\ 
\frac{ A_i^\vee (\bs l+\bs 1_i,k+1)}{\ n\ -\ e_{2-i}^\vee(\bs l+\bs 1_i,k+1)
 }\ b_n (\bs l,k)\ . 
\eean
In particular,
\bea
b_{n} (\bs l,k)\ =\ 
{l_2\choose n}\
\prod_{i=0}^{l_2-n-1}
\prod_{j=1}^{2}\
\frac
{\sum_{s=3-j}^2m_s - l_2 + i + j}
{\sum_{s=3-j}^2m_s - l_2 + i + j + 1 + l_{2-j} - l_{3-j}}\ ,
\eea
where $l_0=k$.

We have
\bea\label{Lin rel}
 \ A_{12}(\bs l,k)\ A_0(\bs l,k)\ y_1(\bs l+\bs 1_0,k+1;x) +
\phantom{aaaaaaaaaaaaaaaaaaaaaaaaaaa}
&&
\notag
\\
 A_{20}(\bs l,k)\ A_1(\bs l,k)\ y_1(\bs l+\bs 1_1,k+1;x)\ +
\phantom{aaaaaaaaaaaaa}
&&
\notag
\\
 A_{01}(\bs l,k)\ A_2(\bs l,k)\ y_1(\bs l+\bs 1_2,k+1;x)\ =\ 0\ ,
%\phantom{aaaaaa}
&&
\notag
\eea
\bea
 A_{12}(\bs l,k)\ A_0^\vee(\bs l,k)\ y_2(\bs l-\bs 1_0,k-1;x) +
\phantom{aaaaaaaaaaaaaaaaaaaaaaa}
&&
\\
 A_{20}(\bs l,k)\ A_1^\vee(\bs l,k)\ y_2(\bs l-\bs 1_1,k-1;x) +
\phantom{aaaaaaaaaaaaaa}
&&
\notag
\\
 A_{01}(\bs l,k)\ A_2^\vee(\bs l,k)\ y_2(\bs l-\bs 1_2,k-1;x)\ =\ 0\ ,
%\phantom{aaaaaa}
&&
\notag
\eea
and
\bea
y_1(\bs l,k;0)\ \neq \ 0\ ,
\qquad
y_1(\bs l,k;1)\ \neq \ 0\ ,
\qquad
y_2(\bs l,k;1)\ \neq \ 0\ .
\eea
%We also have
%\bea
%y_1((l,0),k;x)\ &=&\  P^{-k-1, -m_1-1}_{l}(x)\ ,
%\\
%y_1((l,l),k;x)\ &=&\  P_{l}^{-k-1, -m_1-m_2-2}(x)\ ,
%\eea
%where $P^{\alpha,\beta}_l(x)$ is the monic
%Jacobi polynomial on the interval $[0,1]$.

\bigskip

\noindent
{\it Example.} If $\bs l = (2,1)$, then
\bea
&&
y_1((2,1),k;x) 
= (x-1)^2+
\\
&&
\phantom{aaaa}
 \frac{(k-1)(2k+2m_1+m_2)}{(m_1+k-1)(m_1+m_2+k)}
 (x-1) + \frac{k(k-1)}{(m_2+m_1+k)(m_1+k-1)} ,
\eea
\bea
y_2((2,1),k;x) \ = \
 x\ -\ \frac{m_2(m_1+m_2+1)}{(m_2+2)(k+m_1+m_2)} \ .
\eea
In particular, if $(2m_1 + m_2)^2 + k(4m_1 - m_2^2) = 0$. Then 
\bea
y_1((2,1),k;x)\ = \ y_2((2,1),k;x)^2\ . 
\eea
\hfill
$\square$

\bigskip

For $i = 0, 1, 2$, the operator $D_i(\bs l,k)$ is adjoint to the operator
$-D^\vee_i(\bs l+\bs 1_i,k+1)$ in the following sense:
\bea\label{Adj r=2}
\langle D_i(\bs l,k) f  ,\ g \rangle_{\bs l+\bs 1_i,k+1}\ =\ -\ 
\langle f  ,\ D^\vee_i(\bs l+\bs 1_i,k+1) g \rangle_{\bs l,k}\ ,
\eea
if $f \in V(\bs l,k),$\ and\ $g \in U(\bs l+\bs 1_i,k+1)$.

\bigskip

Let $t_1, \dots , t_{l_1}$ be the roots of the polynomial $y_1(\bs l,k;x)$. 
Let $s_1, \dots , s_{l_2}$ be the roots of the polynomial $y_2(\bs l,k;x)$. 
Assume that $y_1(\bs l,k;x)$ and $y_2(\bs l,k;x)$ do not have multiple roots,
then the roots satisfy the system of equations
\bean\label{BAE r=2}
{}
\\
\frac{m_1}{t_i}\ +\ \frac{k}{t_i-1}\ -\
\sum_{j,\ j\neq i} \frac2{t_i - t_j}\ + \
\sum_{j} \frac 1{t_i - s_j}\ = \ 0\ ,
\qquad 
i = 1, \dots , l_1 \ ,
\notag
\\
\frac{m_2}{s_i}\  -\
\sum_{j,\ j\neq i} \frac2{s_i - s_j}\ + \
\sum_{j} \frac 1{s_i - t_j}\ = \  0 \ ,
\qquad 
i = 1, \dots , l_2 \ .
\notag
\eean
If parameters $\bs m, \bs l, k$ are consistent, then system \Ref{BAE r=1}
has at most one  solution $t_1, \dots , t_{l_1}, 
 s_1, \dots , s_{l_2}$ 
up to permutations of the first and second groups of these numbers.

If parameters $\bs l$ are fixed, then for almost all $\bs m, k$ such that
$\bs m, \bs l, k$ are consistent, the polynomials
 $y_1(\bs l,k;x)$ and $y_2(\bs l,k;x)$ do not have multiple roots.

Formulas \Ref{2' recc} give a recurrent way to generate solutions to
\Ref{BAE r=2}.

\bigskip

Coefficients of the polynomials $y_1(\bs l,k;x)$ and $y_2(\bs l,k;x)$ 
are rational functions of $k, \bs m$.
%which are regular if Re $k < -1$ and Re $m_1+m_2 < -3$. Thus $y_1(\bs l,k;x)$
%and $y_2(\bs l,k;x)$ are well defined for
%these values of parameters.

Assume that Re $k < 0$ and Re $m_1+m_2 < -2$.
For polynomials $f, g$, \
set
\bea
(f,\ g)_{k}\ = \ \int_0^1\ 
 \frac{f(x)\ g(x)}{ x^{m_1+m_2+3}\  (1-x)^{k+1} } \ dx\ .
\eea
For $i = 0,1,2$, the operator $D_i(\bs l,k)$ is adjoint to the operator
$-D^\vee_i(\bs l+\bs 1_i,k+1)$ in the following sense:
\bea
( D_i(\bs l,k) f  ,\ g )_{k+1}\ =\ -\ 
( f  ,\ D^\vee_i(\bs l+\bs 1_i,k+1) g )_{k}\  .
\eea

We have
\bea\label{orth r=2}
\int_0^1\ 
 \frac{y_1(\bs l,k;x) \ x^n}{ x^{m_1+m_2+3}\ (1-x)^{k+1} } \ dx\ 
 =\ 0 \ ,
&\qquad &
n = 1, \dots , l_2\ ,
\\
\int_0^1\ 
 \frac{y_1(\bs l,k;x) \ x^n}{ x^{m_1+2}\ (1-x)^{k+1} } \ dx\ =\ 0\ ,
&\qquad &
n = 1, \dots , l_1-l_2\ .
\eea
We also have
\bea\label{norms r=2}
(y_1(\bs l,k;x), \ y_2(\bs l,k;x))_{k}\ = \ C(\bs l,k)\ ,
\eea
where  
\bea
C(\bs l,k)\ =\  l_2! \
\frac{\Gamma(-k+l_1)\ \Gamma(-m_1 -m_2 -2)}
{\Gamma(-k - m_1 - m_2 - 2 + l_1)}
\ {} \times \hspace{5.2cm}
\\
\prod_{i=0}^{l_2-1}
\frac{(m_2-i)\ (m_1+m_2+1-i)}
{(m_2+ l_1-l_2+1-i) (k+m_1+m_2+2-l_1-i)^2}
\prod_{i=0}^{l_1-l_2-1}
\frac
{m_2+2+i}
{k+m_1+1-l_1-i} .
\eea

\bigskip 

In the next formulas we 
include $\bs m$ in the notation for $y_1(\bs l,k;x)$ and $V(\bs l,k)$.
Then for consistent $\bs m, \bs l, k$, the space
$V(\bs m,\bs l,k)$ has a basis formed by polynomials
\bea
&&
y_1((m_1,m_2),(l_1,l_2),k;x)\ ,
\\
&&
\phantom{aaa} x^{m_1+1}\,y_1((-m_1-2,m_1+m_2+1),(k+l_2-l_1,l_2),k;x)\ ,
\\
&&
\phantom{aaaaaa}
x^{m_1+m_2+2}\,y_1((m_2,-m_1-m_2-3),(k+l_2-l_1,k-l_1),k;x)\ .
\eea

\subsection{Application to the Bethe ansatz in the Gaudin model}
We let $\g = \n_+\oplus\h\oplus\n_-$ be a simple Lie algebra,
\ {}  $\al_1, \dots , \al_r \in \h^*$
simple roots. 

%Let $(x_i)_{i\in I}$ be an orthonormal basis in $\g$, \
%$\Omega =  \sum_{i\in I} x_i\otimes x_i\ \in \g \otimes \g$
%the Casimir element. 

For a $\g$-module $L$ and $\mu \in \h^*$
denote by $L[\mu]$ the weight subspace of $L$ of weight $\mu$ and by
$\Sing\, L[\mu]$ the subspace of singular vectors of weight $\mu$.

Let $n$ be a positive integer and  $\bs \La = (\La_1, \dots , \La_n)$,
$\La_i \in \h^*$, a set of integral dominant weights.
For a
weight $\mu \in \h^*$, let $L_{\mu}$ 
be the irreducible $\g$-module with highest weight $\mu$.
Denote by $L_{\bs \La}$ the tensor product 
$L_{\La_1} \otimes \dots \otimes L_{\La_n}$.

Let $\bs z = (z_1, \dots , z_n)$ be a point in $\C^n$ with distinct coordinates.
Introduce linear operators $K_1(\bs z), \dots , K_n(\bs z)$ on $L_{\bs \La}$ by the formula
\bea
K_i(\bs z)\ = \ \sum_{j,\ j \neq i}\ \frac{\Omega^{(i,j)}}{z_i - z_j}\ , 
\qquad i = 1, \dots , n .
\eea
The operators 
are called {\it the Gaudin Hamiltonians} of the Gaudin model associated with
$L_{\bs \La}$. The  Hamiltonians commute.

The problem is to diagonalize simultaneously
the  Hamiltonians, see \cite{B, BF,  FFR, G, MV1, RV, ScV, V2}.

The  Hamiltonians commute with the $\g$-action on $L_{\bs \La}$.
 It is enough to
diagonalize the  Hamiltonians on the subspaces of singular vectors,\
$\Sing \,L_{\bs \La}[\mu] \subset L_{\bs \La}$. 

The eigenvectors of the Gaudin Hamiltonians are constructed by the Bethe ansatz method.

\bigskip

Fix $\mu = \sum_{s=1}^n\La_s - \sum_{i=1}^r l_i\alpha_i$ where
$l_1, \dots , l_r$ are  non-negative integers. Set
\bea
\bs t\ = \ (t^{(1)}_{1}, \dots , t^{(1)}_{l_1}, t^{(2)}_{1}, \dots , t^{(2)}_{l_2},
\dots , t^{(r)}_{1}, \dots , t^{(r)}_{l_r})\ .
\eea
One defines a suitable rational function $w(\bs t,\bs z)$ with values in $L_{\bs \La}[\mu]$
as in \cite{RV}, cf. \cite{MV2, RSV}.
The function  is symmetric 
with respect to the group $\bs\Sigma_{\bs l} =
\Sigma_{l_1}\times \dots \times \Sigma_{l_r}$ of permutations of coordinates $t^{(i)}_j$ 
with the same upper index. One considers the system of equations
\bean\label{Bethe eqn 1}
-\sum_{s=1}^n \frac{(\Lambda_s, \alpha_i)}{t_j^{(i)}-z_s}\ +\
\sum_{s,\ s\neq i}\sum_{k=1}^{l_s} \frac{(\alpha_s, \alpha_i)}{ t_j^{(i)} -t_k^{(s)}}\ +\
\sum_{s,\ s\neq j}\frac {(\alpha_i, \alpha_i)}{ t_j^{(i)} -t_s^{(i)}}
= 0, 
\eean
where $i = 1, \dots , r$ and $j = 1, \dots , l_i$.
The system is symmetric with respect to   $\bs\Sigma_{\bs l}$.
One shows that if $\bs t^0$ is a solution to \Ref{Bethe eqn 1}, then
$w(\bs t^0,\bs z)$ belongs to $\Sing\,L_{\bs \La}[\mu]$ and $w(\bs t^0,\bs z)$ is an eigenvector of
the Hamiltonians $K_1(\bs z), \dots , K_n(\bs z)$, see \cite{RV}. 

This method of finding eigenvectors is called {\it the Bethe ansatz method}.
System \Ref{Bethe eqn 1} is called {\it the Bethe ansatz equation}, 
the vector $w(\bs t^0,\bs z)$ is called {\it a Bethe vector}.

The Bethe Ansatz Conjecture says that if 
$\La_1, \dots , \La_n$ are integral dominant, the weight
$\mu$ is integral dominant, and 
$z_1, \dots , z_n$ are generic, then the Bethe vectors
form a basis in $\Sing \,L_{\bs \La}[\mu]$. In particular, 
the conjecture implies
that the number of $\bs\Sigma_{\bs l}$-orbits of solutions to \Ref{Bethe eqn 1} 
is not less than the dimension of $\Sing\, L_{\bs \La}[\mu]$.

Assume that $\g = sl_{3}$,\ $n=2$,\ $z_1=0,\ z_2=1$. Assume that for some
non-negative integers $k, m_1, m_2$, the weights 
$\La_1$ and $\La_2$ are such that
\bea
(\La_1, \al_1) = m_1\ ,  
\qquad
(\La_1, \al_2) = m_2\ ,
\qquad
(\La_2, \al_1)  = k\ , 
\qquad
(\La_2, \al_2) = 0 \ .
\eea
Then system \Ref{Bethe eqn 1} becomes system \Ref{BAE r=2}. 

Let $\mu = \La_1+\La_2 - l_1\al_1-l_2\al_2$ be integral dominant, then
the space of highest weight vectors $\Sing \,L_{\La_1}\otimes L_{\La_2}[\mu]$
is one-dimensional, if parameters $\bs m,\bs l, k$ are consistent, and
is zero-dimensional otherwise.

It was shown in \cite{MV1} that if $\mu$ is integral dominant and the numbers
 $t_1, \dots , t_{l_1}$ and $s_1, \dots , s_{l_2}$
form a solution to \Ref{BAE r=2}, then $\bs m,\bs l, k$ are consistent and
\bea
\prod_{i=1}^{l_1}\ (x-t_i) \ = \ y_1(\bs l,k;x)\ ,
\qquad
\prod_{i=1}^{l_2}\ (x-s_i) \ = \ y_2(\bs l,k;x)\ .
\eea
Then we obtain, see the Example in Section \ref{subsection r=2}:

\bigskip

\noindent
{\bf A counterexample to the Bethe Ansatz Conjecture.}
\newline
{\it Let $\bs l = (2,1)$. Let  $\bs m, \bs l, k$ be consistent and
\bean\label{EQ}
(2m_1 + m_2)^2\ +\ k (4m_1 - m_2^2)\ = \ 0\ ,
\eean
then 
$$
y_1((2,1),k;x)\ =\  y_2((2,1),k;x)^2 \ ,
$$ 
and roots of these two polynomials
do not form a solution to  \Ref{BAE r=2}. Hence  system 
\Ref{BAE r=2} does not have solutions
in this case. Hence there is no Bethe vector to generate the one-dimensional space
$\Sing \,(L_{\La_1}\otimes L_{\La_2})[\mu]$.}

For example parameters \ $\bs m= (2,3),\ k = 49, \ \bs l = (2,1)$ \
are consistent and satisfy
\Ref{EQ}.

%\bigskip

%This example also gives a counterexample to one of the main claims, 
%Corollary 3.3, 
%of paper \cite{F} by E. Frenkel.

\section{Spaces of polynomials}\label{spaces sec}
In this section we define our main objects - the spaces of polynomials 
with two singular points $0,1$ and special exponents at $1$.
\subsection{An $sl_{r+1}$ lemma}\label{sl lem}
Fix a natural number $r$ and choose non-negative integers
$m_1,\dots,m_r$.
Choose a non-negative integer $k$.

Consider the complex Lie algebra $sl_{r+1}$. Let $\h$ be its Cartan
subalgebra, let $\al_1, \dots, \al_r \in \h^*$ be its 
simple roots, $(\,,\,)$ the
scalar product on $\h^*$ such that $(\al_i, \al_i) = 2$. For a weight
$\La \in \h^*$ denote by $L_\La$ the irreducible $sl_{r+1}$ module
with highest weight $\La$.

Let the weights
$\La_1,\La_2$ be such that the scalar products with simple roots are given
by 
\bean\label{weights}
 \qquad (\La_1,\al_s)=m_s, \qquad (\La_2,\al_s)=\delta_{s1}k,
\eean
where $s=1,\dots,r$.
Thus $L_{\La_2}$ is the $k$-th symmetric power of the vector
representation and $L_{\La_1}$ can be any irreducible
finite-dimensional $sl_{r+1}$ module.

We describe the decomposition of the tensor product $L_{\La_1}\otimes
L_{\La_2}$.
\begin{lem}\label{decom}
We have an isomorphism of $sl_{r+1}$ modules:
\be
L_{\La_1}\otimes L_{\La_2} = \oplus\ L_{\La_1+\La_2-\sum_{s=1}^{r}
l_s \al_s}
\ ,
\ee
where the sum is over non-negative integers
$l_1,\dots, l_{r}$ such that 
\bean\label{consistent}
k\geq l_1\geq l_2\geq \dots\geq l_r\geq 0, \qquad  l_{s}-l_{s+1}\leq m_{s}
\qquad 
(s=1,\dots,r) .
\eean
% $a_1,\dots,a_{r} \in \Z_{\geq 0}$ such that $a_i\leq m_i$, $i=1,\dots,r$ 
% and $\sum_s a_s \le k$. 
In particular the multiplicity of every module in the right hand side is one. 
%all multiplicities are either $0$ or $1$. 
$\qquad \Box$
\end{lem}

We use the notation $\bs l=(l_1,\dots,l_r)$ and $\bs
m=(m_1,\dots,m_r)$. We call $\bs l$ a {\it partition }
if $l_i$ are integers and $l_1\geq l_2\geq \dots\geq l_r\geq 0$.
We set 
\be
l_0=k, \qquad l_{r+1}=0.
\ee

We call parameters $\bs m,\bs l, k$ {\it consistent} if all 
$m_i,l_i$ and $k$ are non-negative integers satisfying inequalities
\Ref{consistent}.

\subsection{Definition of spaces of polynomials $V(\bs m,\bs l,k)$ and $U(\bs m,\bs l,k)$}\label{space}
Fix a positive number $r$.
Let $V$ be an $(r+1)$-dimensional complex vector space of polynomials 
of one variable $x$. 
Assume that $V$ has no base points,
i.e. for any $z \in \C$ there exists $f \in V$ such that
$f(z) \neq 0$.

The {\it Wronskian polynomial $W(V)$ of $V$} is the Wronskian of any basis
$f_1, \dots , f_{r+1}$ of $V$:
\be
W(V)=W(f_1, \dots , f_{r+1})=\det\left(\frac{d^{i}}{d x^{i}}\ 
f_{j}\right).
\ee
The Wronskian polynomial $W(V)$ is defined up to multiplication by a
non-zero complex number.

A point $z\in \C$ is called a {\it singular point of $V$}, if $z$ is a
root of $W(V)$.

For any $z\in \C$ there exist unique non-negative integers $n_1(z)
, \dots , n_r(z)$ such that for any $f \in V$ the order of zero of $f$
at $z$ belongs to the set
\bean\label{exp at z}
\{\ 0, \ 1 + n_1(z), \ 2 + 
n_1(z) + n_2(z),\ \dots , \ r + \sum_{s=1}^r n_s(z)\ \}.
\eean
A point $z \in \C$ is a singular point of $V$ if and only if
at least one of the numbers $n_1(z), \dots , n_r(z)$ is positive.
The set \Ref{exp at z} is called {\it the set of exponents of $V$ at $z$}.

For $z\in\C$ we define a full flag in $V$
\be
\mc F(z)=\{0=F_0(z)\subset
F_1(z)\subset F_2(z)\subset\dots\subset F_{r+1}(z)=V\},
\ee
where  the $i$-dimensional subspace
$F_i(z)$ consists of all polynomials with root $x=z$ of order at least
$r - i + 1 + \sum_{s=1}^{r-i+1} n_s(z)$. We call the flag $\mc F(z)$ 
{\it the flag of $V$ at $z$}.

Similarly, there exists unique non-negative integers 
$ n_1(\infty),  \dots , n_{r}(\infty), l$ 
such that for any $f \in V$ the degree of $f$ 
belongs to the set
\bean\label{exp at infinity}
\qquad \{\ l,\ l+1+n_1(\infty), \  l+2 + 
n_1(\infty) + n_2(\infty),\ \dots , \ l+r + \sum_{s=1}^r n_s(\infty)\ \}\
.
\eean
The set \Ref{exp at infinity} is called {\it the set of exponents of $V$ at
infinity}. 

We define a full flag in $V$ 
\be
\mc F(\infty)=\{0=F_0(\infty)\subset F_1(\infty)\subset
F_2(\infty)\subset\dots\subset F_{r+1}(\infty)=V\},
\ee
where  the $i$-dimensional subspace
$F_i(\infty)$ consists of all polynomials of degree no greater than
$l+i-1 + \sum_{s=1}^{i-1}n_s(\infty)$. We call the flag $\mc F(\infty)$ {\it the flag of $V$ at infinity}.

We have
\be
W(V)=\prod_{z\in\C}(x-z)^{\sum_{s=1}^r(r+1-s)n_s(z)},
\ee
and therefore the following relation
\bean\label{reso}
\sum_{z\in \C} \sum_{s=1}^r\ (r+1-s)\ n_s(z) \ = \ (r+1)\, l\ + \
\sum_{s=1}^r \ (r+1-s)\ n_s(\infty) \ .
\eean

Given a full flag $\mc F=\{0=F_0\subset F_1
\subset \dots \subset F_{r+1}=V\}$ we call a basis
$\{f_1,\dots,f_{r+1}\}$ of $V$ {\it compatible with flag $\mc F$} if
$\{ f_1,\dots,f_s\}$ is a basis of $F_s$ for $s=1,\dots,r$.

We use the notation $\bs
n(\infty)=(n_1(\infty),\dots,n_r(\infty))$, $\bs n(z)=(n_1(z),\dots,n_r(z))$.

Fix arbitrary sequences of non-negative integers  
$\bs n(z),\bs n(\infty)$ so 
that they
satisfy relation \Ref{reso}.  Let $\La(z)$, $\La(\infty)$ be
$sl_{r+1}$ integral dominant weights defined by the conditions
\be
(\La(z), \al_i) = n_i(z),\qquad (\La(\infty),\al_i) = n_i(\infty)\qquad (i =
1, \dots, r). 
\ee
Then it is known, see \cite{EH}, \cite{MV1}, that the
number of spaces $V$ with the given
$\bs n(z),\bs n(\infty),l$ is finite and bounded from above by the multiplicity
of the module $L_{\La(\infty)}$ in the tensor product
$\otimes_z L_{\La(z)}$. Moreover the number of such spaces is non-zero if
and only if the above multiplicity is non-zero.

\medskip

In this paper we study spaces of polynomials which have only two
singular points $z_1=0$ and $z_2=1$.  Moreover, we consider spaces
only with special exponents at $1$. Namely, we restrict our
attention to two cases
\be
\bs n(1)=(k,0,\dots,0)
\quad {\rm and} \quad
\bs n^\vee(1)=(0,\dots,0,k).
\ee

For our purposes we parameterize the exponents with the above
restriction using the parameters of Section \ref{sl lem}.  
Given parameters $\bs m, \bs l,k$
we define two sets of parameters $\bs n(0),\bs n(1),l,\bs n(\infty)$
and $\bs n^\vee(0),\bs n^\vee(1),l^\vee,\bs n^\vee(\infty)$ by the
formulas:
\begin{align}\label{exp1}
& n_i(1)=\delta_{1i}k, \hspace{20pt} n_i(0)=m_i,  \\
&l=l_1, \hspace{54pt} n_i(\infty)=l_{i-1}+l_{i+1}-2l_i+m_i+\delta_{1i}k\hspace{20pt} 
(i=1,\dots,r) \notag
\end{align}
and
\begin{align}\label{exp2}
&n_i^\vee(1)=\delta_{ir}k,  \hspace{10pt} n_i^\vee(0)=m_{r+1-i},
\\
&l^\vee=l_r,\hspace{30pt}  
n_i^\vee(\infty)=l_{r-i}+l_{r+2-i}-2l_{r+1-i}+m_{r+1-i}+\delta_{ri}k
\hspace{10pt} (i=1,\dots,r). \notag
\end{align}
Both  $\bs n(0),\bs n(1),l,\bs n(\infty)$ and  
$\bs n^\vee(0),\bs n^\vee(1),l^\vee,\bs n^\vee(\infty)$ 
satisfy relation \Ref{reso}.
%It is easy to see that each of the maps $\{\bs m, \bs
%l, k\} \to \{\bs n(0),\bs n(1),l,\bs n(\infty)\}$ and $\{\bs m, \bs
%l, k\} \to \{\bs n^\vee(0),\bs n^\vee(1),l^\vee,\bs n^\vee(\infty)\}$ 
%is bijective (under the assumption of relation \Ref{reso}).

In this parameterization the sets of exponents at infinity given by
\Ref{exp at infinity} take the form $\{e_0(\bs m,\bs
l,k),\dots,e_{r}(\bs m,\bs l,k)\}$ and $\{e_0^\vee(\bs m,\bs
l,k),\dots,e_r^\vee(\bs m,\bs l,k)\}$, where
\begin{align}
%&e_0(\bs m,\bs l,k)=l_1, \quad 
&e_i(\bs m,\bs l,k)=
k+\sum_{s=1}^{i}m_s -l_i+l_{i+1}+i, %\qquad (i=1,\dots,r), 
\label{exp at inf 1} \\
%&e_r^\vee(\bs m,\bs l,k)=k+\sum_{s=1}^rm_s-l_1+r, \quad 
&e_i^\vee(\bs m,\bs l,k)=
\sum_{s=r+1-i}^rm_s-l_{r-i+1}+l_{r-i}+i\label{exp at inf 2},
\end{align}
where $i=0,\dots,r$. Note that for $0\leq i\leq j\leq r$ we have
\bean\label{exp rel}
e_j(\bs m,\bs l,k)-e_i(\bs m,\bs l,k)=
e_{r-i}^\vee(\bs m,\bs l,k)-e_{r-j}^\vee(\bs m,\bs l,k).
\eean
 
Specializing the general statements to our case we obtain:
\begin{lem}\label{V}
  If parameters $\bs m, \bs l, k$ are consistent then there exists a unique
  space of polynomials $V(\bs m,\bs l,k)$ 
  of dimension $r+1$ with exactly two singular
  points $0$ and $1$ and exponents given by \Ref{exp at z} and
  \Ref{exp at infinity} where $\bs n(0),\bs n(1),\bs n(\infty)$
  are given by \Ref{exp1}. If parameters
  $\bs m, \bs l, k$ are not consistent then there exists no such a space. 
  $\qquad \Box$
\end{lem}

\begin{lem}\label{U}
  If parameters $\bs m, \bs l, k$ are consistent then there exists a unique
  space of polynomials $U(\bs m,\bs l,k)$ 
  of dimension $r+1$ with exactly two singular
  points $0$ and $1$ and exponents given by \Ref{exp at z} and
  \Ref{exp at infinity} where $\bs n(0),\bs n(1),\bs n(\infty)$ are
  equal to $\bs n^\vee(0),\bs n^\vee(1),\bs n^\vee(\infty)$ 
  given by \Ref{exp1}. If parameters
  $\bs m, \bs l, k$ are not consistent then there exists no such a space. 
  $\qquad \Box$
\end{lem}

Our parameters $\bs m$ are often fixed. We often drop the dependence 
on $\bs m$ from our notation and simply write $V(\bs l,k)$, $U(\bs l,k)$.

\subsection{Connection between spaces 
$V(\bs m, \bs l,k)$ and $U(\bs m, \bs l,k)$}
Let $\bs m, \bs l, k$ be consistent parameters, 
 $V = V(\bs m, \bs l,k)$ and $U = U(\bs m, \bs l,k)$ 
 the corresponding spaces of polynomials.

We define a tuple of polynomials $\bs T=(T_1,\dots,T_r)$ by 
\bean\label{T}
T_1=(x-1)^kx^{m_1}, \qquad T_i=x^{m_i} \qquad (i=2,\dots,r).
\eean
\begin{lem} We have
\be
W(V)=T_1^rT_2^{r-1}\dots T_r, \qquad W(U)=T_r^rT_{r-1}^{r-1}\dots T_1.
\ee
$\qquad \Box$
\end{lem}

For $f_1, \dots, f_i \in V$ and $g_1,\dots,g_i \in U$  
define the divided Wronskians 
\bea
&&W^\dagger_V(f_1, \dots , f_i)\ =
\ W(f_1,\dots,f_i)\ T_1^{1-i} T_2^{2-i} \dots T_{i-1}^{-1},
\\
&&
W^\dagger_U(g_1, \dots , g_i)\ =\ W(g_1, \dots , g_i)\ 
T_r^{1-i} T_{r-1}^{2-i} \dots T_{r-i+2}^{-1}.
\eea
We also define divided Wronskian $W_V^\dagger (V_1)$ 
(resp. $W_U^\dagger (U_1)$)
of any subspace $V_1\subset V$ (resp. $U_1\subset U$) as 
the divided Wronskian of any basis in $V_1$ (resp. $U_1$). 
The divided Wronskians $W_V^\dagger (V_1)$, $W_U^\dagger (U_1)$ are 
defined up to multiplication by a non-zero constant.

The divided Wronskians are polynomials.

\begin{lem}\label{dual}
We have
\bea
V\ =\ \{W^\dagger_U(g_1,\dots,g_r), \ g_1, \dots , g_r \in U \}\ ,
\\
U\ =\ \{W^\dagger_V(f_1, \dots , f_r), \ f_1, \dots , f_r \in U \}\ .
\eea
\end{lem}
\begin{proof}
  The vector spaces $V$ and $\{W^\dagger_U(g_1, \dots , g_r), \ g_1,
  \dots , g_r \in U \}$ have the same dimensions, the same singular
  points, and exponents, see \cite{MV1} for details. 
  Hence they coincide by Lemma \ref{V}.
\end{proof}

\subsection{Examples of spaces}\label{ex sec}
 Let $m_1,\dots,m_r$ be non-negative numbers.
Let $\bs 0 = (0, \dots , 0)$.
Then the parameters $\bs m, \bs 0, k=0$ are consistent.

Consider the space $V_0 = V(\bs m, \bs 0, 0)$.
Let $e_i=e_i(\bs m, \bs 0, 0)$. Then we have
\be
e_i=\ i + \sum_{j=1}^i m_j,
\ee
where $i = 0, \dots, r$, see \Ref{exp at inf 1}.

\begin{lem}\label{V0} 
We have
\be
V_0={\rm span}\{1=x^{e_0},\
x^{e_1},\ \dots ,\ x^{e_r}\}.
\ee
$\qquad \Box$
\end{lem}

\begin{proof}
  The exponents of the space $V_0$ at $0$ and infinity coincide. Thus
  the polynomial in $V_0$ with the highest order of zero at $0$ has to
  be $x^{e_r}$.  Similarly, a polynomial in $V_0$ with the order
  of zero at $0$ equal to $e_{r-1}$ has to be a linear combination of
  $x^{e_{r-1}}$ and $x^{e_r}$.  Continuing this argument we
  obtain the lemma.
\end{proof}

Fix a non-negative integer $l$, and let 
$\bs l^0 = (l, \dots , l)$. Assume that $m_r\geq l$.
Then the parameters $\bs m, \bs l^0, k=l$ are consistent.

Consider the space $U_l = U(\bs m, \bs l^0, l)$. 
Let $e_i^\vee=e_i^\vee(\bs m, \bs l^0, l)$. Then we have
\be
e_0^\vee=\ l\ , 
\qquad
e_i^\vee=e_i^\vee(\bs m, \bs l^0, l)\ = \ i + \sum_{j=1}^i m_{r+1-j}\ ,
\ee
where $i = 1, \dots, r$, see \Ref{exp at inf 2}.

\begin{lem}
There exists a unique monic polynomial $u_0(l)$ of degree $l$ such that 
\be
U_l={\rm span}\{u_0(l),\  x^{e_1^\vee},\ \dots, \
x^{e_r^\vee}\}.
\ee
Moreover, $u_0(l) \,= \,x^l\, - \,B_l\, x^{l-1} \, +\,
\dots $,  where
\bean\label{Ak}
 B_l\ =\ l\prod_{s=1}^r\frac{e_0^\vee
- e_s^\vee}
{e_0^\vee - e_s^\vee - 1}\ .
\eean
\end{lem}
\begin{proof}
The last $r$ exponents of $U_l$ at 0 are equal to the last
$r$ exponents of $U_l$ at infinity. Thus the functions
$x^{e_1^\vee}$, \dots , $x^{e_r^\vee}$ belong to 
$U_l$. Therefore there exists
a unique monic polynomial $u_0(l)$ of degree $l$ such that the functions
$u_0(l), x^{e_1^\vee}, \dots, x^{e_r^\vee}$ form a basis of
$U_l$.  Moreover, we have
\be
W(u_l, x^{e_1^\vee},  \dots , x^{e_r^\vee})\ =\
c (x-1)^l\ x^{\sum_s e_s^\vee - (r+1)r/2},
\ee
where $c$ is a constant. The first non-trivial coefficient $B_l$ of the 
polynomial $u_0(l)$ is easily computed from this formula.
\end{proof}

\begin{cor} \label{U0}
We have
\be
U_0={\rm span}\{
1 = x^{e_0^\vee(\bs m, \bs 0, 0)},\ 
x^{e_1^\vee(\bs m, \bs 0, 0)},\ \dots ,\ 
x^{e_r^\vee(\bs m, \bs 0, 0)}\}.
\ee
$\qquad \Box$
\end{cor}

\section{Recursion for spaces of polynomials}\label{recursion sec}
In this section we explain how the spaces of polynomials 
with different values of parameters are related to each other.
\subsection{Recursion of the first type}\label{V sec}
Fix a sequence $\bs m$ of non-negative integers.  Let $\bs l,k$ be
such that $\bs m, \bs l, k$ are consistent parameters. 
% Let $V(\bs
%l,k)$ be the corresponding space of polynomials with exponents
%described by \Ref{exp1}, see Lemma \ref{UV}.

Introduce linear first order differential operators
\begin{align}\label{op1}
%&D_0(\bs l,k)=x(x-1)\frac d {dx}-l_1x+k+1, \\ &
D_i(\bs l,k)=x(x-1)\frac d {dx}-(k+\sum_{s=1}^{i}m_s-l_{i}+l_{i+1}+i)(x-1)-k-1,
%\notag 
\end{align}
where $i=0,1,\dots,r$.  Note that the coefficients of $(x-1)$ are exactly
the degrees of polynomials occurring in the space $V(\bs l,k)$.
Namely, the coefficient of $(x-1)$ in $D_i(\bs l,k)$ is equal to 
$e_i(\bs l,k)$, see \Ref{exp at inf 1}.

The operators $D_i(\bs l,k)$ are linearly dependent. Let 
\bean\label{Aab}
A_{ab} = A_{ab}(\bs l,k)=e_a(\bs l, k) - e_b(\bs l, k).
\eean
\begin{lem}
\label{V 2dim}
The vector space spanned by operators $D_i(\bs l,k)$, $i=0,\dots,r$,
has dimension two.
Moreover,  for $0\leq i < j < s\leq r$, we have
\bea
 A_{js} D_i(\bs l,k) \ +\ A_{si} D_j(\bs l,k)\ 
+\ A_{ij} D_s(\bs l,k)\ =\ 0\ .
\eea
$\qquad \Box$
\end{lem}

Set $\bs 1_i=(1,\dots,1,0,\dots,0)$ where we have $i$ ones and $r+1-i$ zeroes.

The operators $D_i$ possess the following holonomy property.
\begin{lem}\label{holon}
For all $i,j\in\{0,1,\dots,r\}$ we have
\be
D_j(\bs l+\bs 1_i,k+1)\ 
D_i(\bs l,k)=D_i(\bs l+\bs 1_j,k+1)\ D_j(\bs l,k).
\ee
\end{lem}
\begin{proof}
It is an easy explicit check.
\end{proof}
We give a conceptual reason for Lemma \ref{holon}
in Section \ref{can map}.

\begin{theorem}\label{V thm}
  Suppose $\bs m, \bs l, k$ and $\bs m, \bs l+\bs 1_i, k+1$ are consistent
  parameters. Then operator $D_i(\bs l,k)$ maps 
$V(\bs l,k)$ to $V(\bs l+\bs 1_i,k+1)$ and
\be
D_i(\bs l,k):\ V(\bs l,k) \to V(\bs l+\bs 1_i,k+1)
\ee
is an isomorphism of vector spaces.
\end{theorem}
\begin{proof}
We have
\begin{align*}
%&\on{Ker}\ D_0 (\bs l,k)= c (x-1)^{k+1}x^{l_1-k-1},\\ &
\on{Ker}\ D_j (\bs l,k)= c(x-1)^{k+1}x^{\sum_{s=1}^{j}m_s-l_{j}+l_{j+1}+j-1},
\end{align*}
where $j=1,\dots,r$ and $c$ is an arbitrary constant. 

Note that the kernel of $D_0$ contains no non-trivial polynomial.
Recall that numbers $\sum_{s=1}^{i-1}m_s+i-1$ and
$\sum_{s=1}^{i}m_s+i$ are consecutive orders of zeros at $x=0$
which appear in $V({\bs l,k})$.  Under our assumptions we have 
\be
\sum_{s=1}^{i-1}m_s+i-1<\sum_{s=1}^{i}m_s-l_{i}+l_{i+1}+i-1
<\sum_{s=1}^{i}m_s+i.  
\ee 
Therefore the intersection of $\on{Ker}\ 
D_i(\bs l,k)$ and $V(\bs l,k)$ is trivial, because the orders of zero
at $x=0$ of functions in these two spaces do not match.

It follows that $\bar V:=D_i(\bs l,k)(V(\bs l,k))$ is a space of polynomials of
dimension $r+1$. Let us consider the exponents of $\bar V$.

Let $f_1,\dots,f_{r+1}$ be a basis in $V(\bs l,k)$ compatible with the
flag at infinity. We have $\deg f_s=e_s$, where $e_s$ are given 
by \Ref{exp at inf 1}. Then we obviously have
\begin{align*}
&\deg D_i f_s=\deg f_s +1 \qquad (s\neq i), \\
&\deg D_i f_i\leq \deg f_i. 
\end{align*}
Also if $g\in V$ and $g(1)=0$ then the order of vanishing at $x=1$ of
$D_i(\bs l,k)g(0)$ is at least $k+2$.

It follows from relation \Ref{reso} that the above inequality is an 
equality and 
the exponents of $\bar V$
coincide with the exponents of $V(\bs l+\bs 1_i,k+1)$. Therefore those
two spaces coincide by the uniqueness part of Lemma \ref{V}.
\end{proof}

\begin{cor}\label{ph-nonph} Let $\bs m, \bs l, k$ be consistent and 
$\bs m,\bs l+\bs 1_i,k+1$ are not consistent. Then $l_i-l_{i+1}=m_i$ and
\be
(x-1)^{k+1}x^{\sum_{s=1}^{i}m_s-l_{i}+l_{i+1}+i-1}\ \in \ V(\bs l,k). 
\ee
In particular $\dim D_i(\bs l,k) V(\bs l,k)=r$.
\end{cor}
\begin{proof}
  If the statement of the lemma were wrong, then by the same argument
  as in the proof of Theorem \ref{V thm}, the space $D_i(\bs l,k) V(\bs l,k)$
  would be a space of polynomials $V(\bs l+\bs 1_i,k+1)$, which does
  not exist by Lemma \ref{V}.
\end{proof}

\begin{lem}\label{comb} 
  Let $\bs m, \bs l, k$ be consistent
  parameters.
\begin{itemize}
 \item There exist a sequence of indices
  $i(1),\dots,i(k)\in \{0,\dots,r\}$ such that 
  $\bs l= \sum_{s=1}^k \bs 1_{i(s)}$ and
 $\bs m, \sum_{s=1}^j
  \bs 1_{i(s)},j$ are consistent parameters for $j=0,\dots,k$. 

 \item This sequence is unique up to 
  permutation of its elements.

 \item $0$ occurs $k-l_1$ times and
  $i$ occurs $l_i-l_{i+1}$ times  $(i=1,\dots,r)$.
\end{itemize}
\end{lem}
\begin{proof}
  Note that the addition of $\bs 1_i$ adds one to $l_i-l_{i+1}$ and
  leaves all other differences $l_j-l_{j+1}$, $j\neq i$ intact. In
  particular the addition of $\bs 1_0$ does not change either of these
  differences. The lemma follows.
\end{proof}

\begin{cor}\label{Vk op}
 Let $\bs l, k$ be such that $\bs m, \bs l, k$ are consistent. There exist 
a unique linear differential operator $D_{\bs l,k}$ of order $k$ and 
of the form 
\begin{align}\label{DD}
D_{\bs l,k}:= D_{j(k)}(\sum_{s=1}^{k-1}\bs 1_{j(s)},k-1)\ \dots\ 
D_{j(2)}(\bs 1_{j(1)},1)\ D_{j(1)}(\bs 0,0)
\end{align}
with $j(s)\in \{0,\dots,r\}$, $\sum_{s=1}^{k}\bs 1_{j(s)}=\bs l$, 
such that $V(\bs l,k)= D_{\bs l,k}\ V(\bs 0,0)$.
\end{cor}
\begin{proof}
  The existence follows immediately from Theorem \ref{V thm} and Lemma
  \ref{comb}. The uniqueness follows from Lemma
  \ref{comb} and Corollary \ref{ph-nonph}.
\end{proof}

Sometimes we stress the dependence of the operator $D_{\bs l,k}$ on
$\bs m$ and write $D_{\bs m,\bs l,k}$.

\medskip

We give some remarks about a basis in $V(\bs l,k)$.
\begin{lem}\label{basis}  The space $V( \bs l,k)$ has a basis
$\{\tilde v_0,\tilde v_1, ...,\tilde v_r\}$ such that $v_0$ 
is a polynomial of degree $l_1$ and
$\tilde v_i$, $i=1,\dots,r$,
have the form
\be
\tilde v_i\ =\ (x-1)^{k+1}\
x^{\sum_{j=1}^{i-1} m_{j}+i-1}\
\sum_{j=0}^{m_i-l_i+l_{i+1}}\
a_{ij}\ x^j,
\ee
where $a_{ij}$ are non-zero constants.
\end{lem}
\begin{proof}
The lemma is true for $V_0$,
\be
x^{e_i}-x^{e_{i-1}}\ =\
(x-1)\ x^{\sum_{j=1}^{i-1} m_{j}+i-1}\
\sum_{j=0}^{m_i}x^j\ .
\ee
In addition if $D = x(x-1)d/d x + A(x-1) -(k+1)$, then
\bean\label{simple act}
D\ :\ (x-1)^{k+1}x^a \ \mapsto\
(k+1+a+A)\ (x-1)^{k+2}x^a \ .
\eean
These two statements imply the lemma.
\end{proof}
Note that using \Ref{simple act}, one can write
the constants $a_{ij}$ explicitly 
in a factored form.

Note also that Corollary \ref{ph-nonph} follows from Lemma \ref{basis}.

%----------------------------

\subsection{Recursion of the second type}\label{U sec}
The construction of this section is parallel to that of Section \ref{V sec}.

Fix a sequence $\bs m$ of non-negative integers.
Let $\bs l,k$ be such that $\bs m, \bs l, k$ are consistent parameters.
% and let
%$U(\bs l,k)$ be the corresponding space of polynomials  
%with exponents described by \Ref{exp2}, see Lemma \ref{UV}.
Introduce linear first order differential operators
\begin{align}\label{U op}
%  &D^\vee_0(\bs l,k)=x\frac d {dx}-(k+\sum_{s=1}^{r}m_s-l_{1}+r+1) , \\ &
  D^\vee_i(\bs l,k)=
  x\frac d {dx}-(\sum_{s=i+1}^{r}m_s-l_{i+1}+l_{i}+r-i), %\notag
\end{align}
where $i=0,\dots,r$.  Note that the constant coefficients in operators
$D^\vee_i(\bs l,k)$ are exactly the degrees of polynomials occurring
in the space $U(\bs l,k)$. Namely, the constant coefficient
in $D^\vee_i(\bs l,k)$ is exactly $e^\vee_{r-i}(\bs
l,k)$, see \Ref{exp at inf 2}.

The operators $D_i^\vee$ are linearly dependent. Recall constants $A_{ab}$ given by \Ref{Aab}.
\begin{lem}
\label{U 2dim}
The vector space spanned by operators $D_i^\vee(\bs l,k)$, $i=0,\dots,r$, 
has dimension two.
Moreover,  for $0\leq i < j < s\leq r$, we have
\bea
 A_{js} D_i^\vee(\bs l,k) \ +\ A_{si} 
D_j^\vee(\bs l,k)\ +\ A_{ij} D_s^\vee(\bs l,k)\ 
=\ 0\ .
\eea
$\qquad \Box$
\end{lem}

The operators $D_i^\vee$ possess the following holonomy property.
\begin{lem}\label{U holon}
For all $i,j\in\{0,1,\dots,r\}$ we have
\be
D_j^\vee(\bs l-\bs 1_i,k-1)\ D_i^\vee(\bs l,k)=
D_i^\vee(\bs l-\bs 1_j,k-1)
\ D_j^\vee(\bs l,k).
\ee
$\qquad \Box$
\end{lem}
We give a conceptual reason for Lemma \ref{U holon}
in Section \ref{can map}.

\begin{theorem}\label{U thm}
  Suppose $\bs m, \bs l-\bs 1_i, k-1$ and $\bs m, \bs l,k$ are consistent
  parameters. Then operator  $D_i^\vee(\bs l,k)$ maps 
  $U(\bs l,k)$ to $U(\bs l-\bs 1_i,k-1)$ and 
\be
D_i^\vee(\bs l,k):\ U(\bs l,k) \to U(\bs l-\bs 1_i,k-1)
\ee
is an isomorphism of vector spaces. 
\end{theorem}
 \begin{proof}
The proof is parallel to the proof of Theorem \ref{V thm}.
\end{proof}

\begin{cor}
  Let $\bs l, k$ be such that parameters $\bs m,\bs l, k$ are
  consistent.  There exist a unique linear differential
  operator $D_{\bs l,k}^\vee$ of order $k$ and of the form
\begin{align*}
D_{\bs l,k}^\vee:=D_{j(1)}^\vee(\bs 1_{j(1)},1)\ 
D_{j(2)}^\vee(\bs 1_{j(1)}+\bs 1_{j(2)},2)
\ \dots\  D_{j(k)}^\vee(\sum_{s=1}^{k}\bs 1_{j(s)},k)
\end{align*}
with $j(s)\in\{0,\dots,r\}$, $\sum_{s=1}^{k}\bs 1_{j(s)}=\bs l$,
such that $U(\bs 0,0)= D_{\bs l,k}^\vee\ U(\bs l,k) $. $\qquad \Box$
\end{cor}

\subsection{Canonical maps}\label{can map}
Fix a sequence $\bs m$ of non-zero integers.  Let $\bs l_1,k_1$ and $
\bs l_2,k_2$ be such that the parameters $\bs m,\bs l_1,k_1$ and $\bs
m, \bs l_2,k_2$ are consistent.

Recall that given any space of polynomials $V$ and $z\in\C$ we have the full
flag of $V$ at $z$, 
$\mc F(V,z)=\{0\subset F_1(V,z)\subset F_2(V,z)\subset \dots
\subset F_{r+1}(V,z)=V\}$, defined by the order of vanishing of
polynomials in $V$ at $x=z$. We also have the full flag of $V$ at infinity,
$\mc F(V,\infty)=\{0\subset F_1(V,\infty)\subset F_2(V,\infty)\subset \dots
\subset F_{r+1}(V,\infty)=V\}$, defined by the degrees of polynomials
in $V$.

We call a linear map $a:\ V(\bs l_1,k_1)\to V(\bs l_2,k_2)$
{\it canonical} if $a$ is an isomorphism and it is compatible with flags at
singular points and infinity in the following sense: 
\begin{align*}
&a(F_i(V(\bs l_1,k_1),0))=F_i(V(\bs l_2,k_2),0) \hspace{35pt}  
(i=1,\dots,r),\\
&a(F_i(V(\bs l_1,k_1),\infty))=F_i(V(\bs l_2,k_2),\infty)\qquad   
(i=1,\dots,r),\\
&a(F_r(V(\bs l_1,k_1),1))=F_r(V(\bs l_2,k_2),1).
\end{align*}

We call a linear map $b:\ U(\bs l_1,k_1)\to U(\bs l_2,k_2)$
{\it canonical} if $b$ is an isomorphism and it is compatible with flags at
singular points and infinity in the following sense: 
\begin{align*}
&b(F_i(U(\bs l_1,k_1),0))=F_i(U(\bs l_2,k_2),0) \hspace{35pt}  
(i=1,\dots,r),\\
&b(F_i(U(\bs l_1,k_1),\infty))=F_i(U(\bs l_2,k_2),\infty)\qquad   
(i=1,\dots,r),\\
&b(F_1(U(\bs l_1,k_1),1))=F_1(U(\bs l_2,k_2),1).
\end{align*}

Clearly inverse of a canonical map is canonical and a composition of
canonical maps is canonical. A canonical map multiplied by a non-zero
number is canonical.

\begin{theorem}\label{can thm}
  Let $\bs m,\bs l_1,k_1$ and $\bs m, \bs l_2,k_2$ be consistent
  parameters. Then there exist unique (up to multiplication by a
  non-zero number) canonical maps $V(\bs l_1,k_1)\to V(\bs l_2,k_2)$
  and $U(\bs l_1,k_1)\to U(\bs l_2,k_2)$.
\end{theorem}
\begin{proof}
  We give a proof for the case of spaces $V(\bs l,k)$. The spaces
  $U(\bs l, k)$ are treated similarly.
  It is obviously sufficient 
  to consider the case of $(\bs l_1,k_1)=(\bs 0,0)$, that is the case of
  $V(\bs l_1,k_1)=V_0$.
  
  We have a map $V_0\to V(\bs l_2,k_2)$ given by the operator
  $D_{\bs l_2,k_2}$, see Corollary \ref{Vk op}. This map is clearly
  canonical, see the proof of Theorem \ref{V thm}.
  
  Since we have the existence, to prove that the canonical map is
  unique it is sufficient to consider canonical maps $V_0 \to V_0$. 
  
  Let $a:\ V_0\to V_0$ be a canonical map. Let $\mc F(0), \mc F(1),\mc
  F(\infty)$ be the corresponding flags at $0,1$ and $\infty$. By
  degree reasons we have
\bean\label{1}
W^\dagger_{V_0} (F_i(\infty))=1 \qquad (i=1,\dots,r)
\eean
(up to non-zero constants).

For $i=1,\dots,r+1$ consider $G_i=F_i(\infty)\cap F_{r+2-i}(0)$.  By the
definition of a canonical map,  $a(G_i)=G_i$.  But because of \Ref{1} we
have $\dim G_i=1$. Let $G_i=cf_i$, $f_i\in V_0$. The basis 
$\{f_1,\dots,f_{r+1}\}$ is compatible with $\mc F(\infty)$ and the basis 
$\{f_{r+1},\dots,f_{1}\}$ is compatible with  $\mc F(0)$.
In the basis
$\{f_1,\dots,f_{r+1}\}$ the map $a$ is diagonal. Let $a (f_s)=\mu_s f_s$.
  
Similarly for $i=1,\dots,r+1$ we have $F_i(\infty)\cap
F_{r+2-i}(1)=cg_i$ for some $g_i\in V_0$. The basis 
$\{g_1,\dots,g_{r+1}\}$ is compatible with $\mc F(\infty)$ and the basis 
$\{g_{r+1},\dots,g_{1}\}$ is compatible with  $\mc F(1)$. The span of 
$g_2,\dots,g_s$ is the intersection of $F_r(1)$ and $F_s(\infty)$ and 
therefore it is preserved by the map $a$.

Without loss of generality we can put $f_1=g_1=1$.  
It follows from Lemma \ref{V0} that without loss of generality we can assume 
\be
f_i=x^{\sum_{s=1}^{i-1}m_s+i-1}.
\ee
It follows that 
if we write $f_i=\sum_{s=1}^i \nu_{si} g_s$ then $\nu_{ii}\neq 0$ and 
$\nu_{1i}\neq 0$.
 
Since the map $a$ preserves the span of $g_2$ we get
$\mu_1=\mu_2=\mu$.  Since the map $a$ preserves the span of $g_2,g_3$
and this span does not contain $g_1$, we get $\mu_3=\mu$. And so on.

It follows that the the map $a$ is scalar.
\end{proof}
Note that holonomy properties described in Lemmas \ref{holon} and
\ref{U holon} are natural to expect because of Theorem \ref{can thm},
since both compositions are the unique canonical maps when acting on
our spaces of polynomials.

\subsection{Inclusion Property} Fix a sequence $\bs m$
of $r$ non-negative integers.  Let a sequence $\bs l$ of $r$
non-negative integers and a non-negative integer $k$ be such
that $\bs m,\bs l, k$ are consistent parameters. Let $s$ be a natural number.
Let $\bar {\bs
  l}=(l_1,\dots,l_r,0,\dots,0)$ where we have $s$ zeroes, and $\bar
{\bs m}=(m_1,\dots,m_r,n_1,\dots,n_s)$, where $n_i$ are non-negative
integers.

\begin{lem}\label{V incl} We have
\be
V(\bs m,\bs l, k)\subset V(\bar{\bs m},\bar{\bs l}, k).
\ee
\end{lem}
\begin{proof}
The lemma is true for the case of $\bs l=0$, $k=0$, see Lemma
\ref{V0}. The canonical isomorphisms
$a:V(\bs m,\bs 0,0)\to V(\bs m, \bs l,k)$ and  
$\bar a:V(\bar{\bs m},\bar{\bs 0},0)\to V(\bar{\bs m}, \bar{\bs l},k)$
are given by the same formula
$a=D_{\bs m, \bs l, k}=D_{\bar{\bs m}, \bar{\bs l}, k}=\bar a$. 
The lemma follows.
\end{proof}

Canonical maps are compatible with the inclusion described in 
Lemma \ref{V incl}. 

Namely, let $\bs l^1, k^1$ be such that $\bs m,\bs l^1,k^1$ are also 
consistent parameters. We set
$\bar {\bs l}^1=(l_1^1,\dots,l_r^1,0,\dots,0)$ 
where we have $s$ zeroes. We have canonical isomorphisms
\bea
a: V(\bs m,\bs l,k)\to V(\bs m,\bs l^1,k^1), 
\qquad %b: U(\bs m,\bs l,k)\to U(\bs m,\bs l^1,k),\\ 
\bar a: V(\bar{\bs m},\bar{\bs l},k)\to V(\bar{\bs m},\bar{\bs l}^1,k^1), 
\qquad 
%\bar b: U(\bar{\bs m},\bar{\bs l},k)\to U(\bar{\bs m},\bar{\bs l}^1,k).
\eea
defined by 
\bea
a=D_{\bs m, \bs l^1,k^1}\circ (D_{\bs m,\bs l,k})^{-1}, \qquad 
%b=(D_{\bs m, \bs l^1,k}^\vee)^{-1}\circ D_{\bs m, \bs l,k}^\vee, \\
\bar a=
D_{\bar{\bs m},\bar{\bs l}^1,k^1}\circ (D_{\bar{\bs m},\bar{\bs l},k})^{-1}. 
%\qquad 
%\bar b=
%(D_{\bar{\bs m},\bar{\bs l}^1,k}^\vee)^{-1}\circ 
%D_{\bar{\bs m},\bar{\bs l},k}^\vee.
\eea

\begin{cor}
The canonical map $a$ is the
restriction  of the canonical map
$\bar a$ to
$V(\bs m,\bs l,k)\subset V(\bar{\bs m},\bar{\bs l},k)$.
$\qquad \Box$
\end{cor}

Spaces $U(\bs m,\bs l,k)$ possess similar properties. 
For example, the analog of Lemma \ref{V incl} 
for spaces $U(\bs m,\bs l,k)$ has the following form. 
\begin{lem}\label{U incl}
We have the surjective well-defined map of vector spaces of polynomials
\bea
\pi: \ U(\bar{\bs m},\bar{\bs l}, k)&\to& U(\bs m,\bs l, k),\\
 g &\mapsto&  x^{-n_1}\frac{d}{d x}\ \ ...\ x^{-n_{s-1}}
\frac{d}{d x}\ x^{-n_s}\frac{d}{d x}\ g.
\eea
$\qquad\Box$
\end{lem}

\subsection{Duality}\label{dual sec}
Fix a sequence $\bs m$ of non-negative integers. Let $\bs l,k$ be such that
$\bs m,\bs l, k$ are consistent parameters. 

We use Lemma \ref{dual} to 
define a pairing $\langle\ ,\ \rangle_{\bs l,k}: V(\bs l,k)\otimes U (\bs l,k)\to \C$ as follows. 
For $f\in V(\bs l,k)$, $g\in U(\bs l,k)$ we set
\be
\langle f, g\rangle_{\bs l,k}=W_{V(\bs l,k)}^\dagger(f, f_1,\dots,f_r), 
\qquad {\rm if} \ \ g=W_{V(\bs l,k)}^\dagger(f_1,\dots,f_r),
\ee
see Section 6.1 in \cite{MV1}, 
where this pairing is studied in a more general setting.
\begin{lem}\label{pair}
The pairing 
$\langle\ ,\ \rangle_{\bs l,k}: V(\bs l,k)\otimes U(\bs l,k)\to \C$ 
is a well-defined 
bilinear, non-degenerate pairing.
\end{lem}
\begin{proof}
Lemma \ref{pair} follows from Lemma \ref{dual}.
\end{proof}
Recall that the degrees of polynomials $e_i$ and $e_j^\vee$ occurring 
in $V(\bs l,k)$ and $U(\bs l,k)$ 
are given by \Ref{exp at inf 1} and \Ref{exp at inf 2}.  
\begin{lem}\label{dot form}
Let $f,g$ be monic polynomials  such that $f\in V(\bs l,k)$, $\deg f=e_i$, and 
$g\in U(\bs l,k)$, $\deg g=e_{r-j}^\vee$. Then we have
\be
\langle f, g\rangle_{\bs l,k}=\prod_{s=0,\ s\neq i}^r(e_s-e_i).
\ee
\end{lem}
\begin{proof}
For $a_1\leq \dots\leq a_j$ we have the formula 
\be
  W(x^{a_1},\dots,x^{a_j})=x^{\sum_{s=1}^ja_s -j(j+1)/2}\prod_{s>t}(a_s-a_t).
\ee
We compute (denoting the terms of lower degrees by dots):
\begin{align*}
&\langle f, g\rangle_{\bs l,k}= \\
&W_{V(\bs l,k)}^\dagger (f, x^{e_0}+\dots,\ 
\dots\ ,\widehat{x^{e_{r-j}}+\dots},\ \dots\ ,x^{e_r}+\dots)
\left(\prod_{s>t, \ s\neq r-j, \atop t\neq r-j}( e_s-e_t)\right)^{-1}= \\
&= (-1)^i \prod_{s>t} (e_s-e_t)\left(\prod_{s>t, 
\ s\neq i, \ t\neq i}( e_s-e_t)\right)^{-1}=
\prod_{s=0,\ s\neq i}^r(e_s-e_i).
\end{align*}

\end{proof}

\begin{prop}
The linear operators $D_i(\bs l,k):\ V(\bs l,k) \to V(\bs l+\bs 1_i,k+1)$ and 
$-D_i^\vee(\bs l+\bs 1_i,k+1):\ U(\bs l+\bs 1_i,k+1)\to U(\bs l,k)$ 
are adjoint. Namely for any $f\in V(\bs l,k)$ and any $g\in U(\bs l,k)$ we have
\be
\langle D_i(\bs l,k) f,g \rangle_{\bs l+\bs 1_i,k+1}= -\langle
 f,D_i^\vee(\bs l+\bs 1_i,k+1)g\rangle_{\bs l,k}.
\ee
\end{prop}
\begin{proof}
  The operator $D_i(\bs l,k)$ defines a canonical map.  It follows
  from definitions that the map adjoint to a canonical map is
  canonical. Therefore, by the uniqueness of the canonical map, the
  map defined by $D_i^\vee(\bs l+\bs 1_i,k+1)$ is a constant multiple 
  of the map
  adjoint to the one defined by $D_i(\bs l,k)$. 
  
  We claim the constant is $-1$. Let $f_0,\dots,f_{r}$ be a basis 
  of monic polynomials in
  $V(\bs l,k)$ compatible with the flag at infinity. 
  We have $\deg f_s=e_s$ and $f_i$ are monic.
  Let $\tilde e_s =e_s+1-\delta_{i,s}$ be the 
  degrees of polynomials in $V(\bs l+\bs 1_i,k+1)$. Let $g\in
  U(\bs l+\bs 1_i,k+1)$ be a monic polynomial of maximal degree.  
  In the next calculation we assume that $i>0$,
  the case of $i=0$ is similar. By Lemma \ref{dot form},
  \begin{align*}
\langle D_i(\bs l,k) f_0,g \rangle&=(e_0-e_{i})
\langle x^{l_1+1}+\dots,g\rangle= \\
&=(e_0-e_{i})\prod_{s,\ s\neq 0}(\tilde e_s-\tilde e_0)=
(e_0-e_{i}+1)\prod_{s=1}^{r}(e_s-e_0).
\end{align*}
Here we denoted  the terms of lower degrees by dots. Similarly we
obtain 
\be 
\langle  f_0, D_i^\vee(\bs l+\bs 1_i,k+1)g
\rangle=(e_{r}^\vee-e^\vee_{r-i}-1)\prod_{s=1}^{r}(e_s-e_0). 
\ee
Therefore by \Ref{exp rel} we have $\langle D_i(\bs l,k) f_0,g \rangle=-\langle
 f_0,D_i^\vee(\bs l+\bs 1_i,k+1)g\rangle$ and the constant is $-1$.
\end{proof}

%---------------------------------------

\section{Spaces $V(\bs m,\bs l,k)$ and Jacobi-Pi\~neiro polynomials}
\label{V sect}
In this section we present an explicit basis of spaces $V(\bs m,\bs l,k)$
in terms of multiple orthogonal polynomials, 
called Jacobi-Pi\~neiro polynomials.
\subsection{Rodrigues formula}
Let $k,l_1,\dots,l_r$ be non-negative integers such that $k\geq
l_1\geq \dots \geq l_r$. Let $m_1,\dots,m_r$ be generic complex
numbers.

For $i=0,\dots,r$ we define a differential operator 
which acts on functions of $x$:
\be
\tilde D_i(\bs m,\bs l,k)=x^{\sum_{j=1}^i{m_j}+i }
\frac{d^{l_i-l_{i+1}}}
{d x^{l_i-l_{i+1}}}x^{l_i-l_{i+1}-\sum_{j=1}^i{m_j}-i}.
\ee
This operator is the composition of multiplication by a monomial,
multiple differentiation, and another multiplication by a monomial.

\begin{lem}\label{commute}  The differential
operators $\tilde D_i(\bs m,\bs l,k)$ commute
for all values of parameters:
\be
\tilde D_i(\bs m_1,\bs l_1,k_1)\tilde D_j(\bs m_2,\bs l_2,k_2)=
\tilde D_j(\bs m_2,\bs l_2,k_2)\tilde D_i(\bs m_1,\bs l_1,k_1).
\ee
\qquad $\Box$
\end{lem}

Set
\bean\label{weight}
\om(x):=(x-1)^{-k-1}x^{-\sum_{i=1}^rm_i-r}.
\eean

Recall operator $D_{\bs m,\bs l,k}$ given by \Ref{DD} and 
our conventions $l_0=k$ and $l_{r+1}=0$.
\begin{prop}\label{Rod}
  We have the equality of differential operators
 \bea 
&&D_{\bs m,\bs l,k}= (x-1)^{k+1}\tilde D_r(\bs m,\bs l,k)\tilde D_{r-1}(\bs m,\bs l,k)\dots \tilde D_0(\bs m,\bs l,k)\frac{1}{x-1}= \\
&&\om^{-1}(x)\ 
  \frac{d^{l_r-l_{r+1}}}{d x^{l_r-l_{r+1}}}\ 
  x^{l_r-l_{r+1}-m_r-1}\frac{d^{l_{r-1}-l_r}}{d
    x^{l_{r+1}-l_r}}\ x^{l_{r-1}-l_{r}-m_{r-1}-1}\dots\ \\ &&\hspace{7.5cm}
  x^{l_1-l_2-m_1-1}\frac{d^{l_0-l_1}}{d x^{l_{0}-l_1}}\ 
  \frac{x^{l_0-l_1}} {x-1}\ .  
\eea
\end{prop}
\begin{proof}
For any complex $a,b$ and non-negative integer $s$ 
we have the following identity of differential operators: 
\begin{align*}
&\left(x(x-1)\frac{d}{d x}-a(x-1)-b-s\right)
\left(x(x-1)\frac{d}{d x}-a(x-1)-b-s+1\right)\dots\\
&{}\hspace{5pt} \dots\left(x(x-1)\frac{d}{d x}-a(x-1)-b\right)=
(x-1)^{b+s+1}x^{a-b+1}\ \frac{d^{s+1}}{d x^{s+1}}
(x-1)^{-b}x^{-a+b+s}.
\end{align*}
Indeed, both operators are of order $s+1$ with the symbol $(x(x-1))^{s+1}d^{s+1}$. Moreover, it is easy to check that both operators kill the functions
\be
(x-1)^{b+i}x^{a-b-s} \qquad (i=0,\dots,s).
\ee
It follows that these two operators are the same.

The proposition follows from the above identity.
\end{proof}

\subsection{Weyl group action}
Let $k,l_1,\dots,l_r$ be non-negative integers satisfying the conditions $k\geq
l_1\geq \dots \geq l_r$. Let $m_1,\dots,m_r$ be generic complex
numbers.  

Let $\mc S_{r+1}$ be the the Weyl group of $sl_{r+1}$.
Then $\mc S_{r+1}$ is isomorphic to the
group of permutations of the set $\{0,1\dots,r\}$.
We denote $w(i)$ the image of $i\in \{0,1\dots,r\}$ under permutation $w$.
The group $\mc S_{r+1}$ is 
generated by simple transpositions $s_1,\dots,s_r$. 
The simple transposition $s_i$ permutes elements $i$ and $i-1$ and 
preserves all other elements of the set $\{0,1,\dots,r\}$.

The group $\mc S_{r+1}$ acts on $sl_{r+1}$ weights by the formula
\be
s_i\cdot \la=\la-(\la+\rho,\al_i)\al_i.
\ee
The first formula in \Ref{weights} identifies  $sl_{r+1}$ weights
with vectors $\bs m$. Under this identification the Weyl group action
becomes:
\be
\s_i\cdot \bs m=
(m_1,\dots,m_{i-2}, m_{i-1}+m_i+1,-m_i-2,m_{i+1}+m_i+1,m_{i+2},\dots,m_r).
\ee
We denote the result of the action of $w\in\mc S_{r+1}$ on 
$\bs m$ by $w\cdot\bs m$.

Given $k$, define another action of the Weyl group
$\mc S_{r+1}$ on the sequences of numbers
$\bs l$
by the formulas
\be
(s_i)_k\bs l=(l_1,\dots,l_{i-1}, l_{i-1}+l_{i+1}-l_i,l_{i+1},\dots,l_r),
\ee
where as always $l_0=k$ and $l_{r+1}=0$.
We denote the result of the 
action of $w\in\mc S_{r+1}$ on $\bs l$ by $(w)_k\bs l$.

If we consider the map $\C^r \to \C^{r+1}$ mapping $\bs l \to 
\{l_0-l_1, l_1-l_2,\dots,l_{r}-l_{r+1}\}$, 
then the second Weyl group action descends
to the standard permutation action of coordinates.

The Weyl group acts on operators $\tilde D_i(\bs m, \bs l,k)$ 
as follows.
\begin{lem}\label{Weyl and D}
We have
\bea
\tilde D_i(s_1\cdot \bs m, (s_1)_k\bs l,k)&=&x^{-m_1-1}\tilde 
D_{s_1(i)}(\bs m, \bs l, k)\ x^{m_1+1},\\
\tilde D_i(s_j\cdot\bs m, (s_j)_k\bs l,k)&=&\tilde D_{s_j(i)}(\bs m, \bs l, k)
,
\eea
where $j=2,\dots,r$. $\qquad \Box$
\end{lem}

\subsection{Basis in $V(\bs m,\bs l,k)$}
Let $\bs m, \bs l, k$ be consistent parameters.

The following lemma follows from Theorem \ref{V thm} and Lemma \ref{V0}.
\begin{lem}\label{basis lem}
The space  $V(\bs m,\bs l,k)$ has a basis of the form
\be 
\{v_0(\bs m, \bs l,k),\ v_1(\bs m, \bs l,k)x^{m_1+1},\dots, v_r(\bs m, \bs l,k)x^{\sum_{i=1}^rm_i+r}\},
\ee
where $v_i(\bs m, \bs l,k)$ is a monic polynomial of degree
\be
\deg v_i(\bs m, \bs l,k)=k-l_i+l_{i+1}
\ee
whose coefficients are rational functions of $\bs m$ 
with coefficients in $\Q$.
Moreover, the polynomial $v_i(\bs m, \bs l,k)$ is obtained  
via the action of operator $D_{\bs m,\bs l,k}$ as follows
\bean\label{ui}
v_i(\bs m, \bs l,k)=
x^{-\sum_{j=1}^im_j-i}\ D_{\bs m,\bs l,k}\  
\left( x^{\sum_{j=1}^im_j+i}\right).
\eean
$\qquad \Box$
\end{lem}

The polynomials $v_i$ are permuted by the Weyl group.
\begin{prop}\label{v inv}
Let $w\in \mc S_{r+1}$ be an element of the Weyl group, then
\be
v_i(w\cdot \bs m, (w)_k \bs l,k)=v_{w(i)}(\bs m,\bs l,k).
\ee
\end{prop}
\begin{proof}
For $i=2,\dots,r$  %we have 
%\be
%s_i \cdot (\sum_{t=1}^jm_t+j)=\sum_{t=1}^{s_i(j)}m_t+s_i(j)
%\ee
%and 
from Lemmas \ref{Weyl and D} and \ref{commute} we have 
the equality of operators
\be
D_{s_i\cdot \bs m,(s_i)_k\bs l,k}=D_{\bs m,\bs l,k}
\ee
%Similarly, we have 
%\be
%s_1 \cdot (\sum_{t=1}^jm_t+j)=-m_1-1+\sum_{t=1}^{s_1(j)}m_t+s_1(j)
%\ee
and 
\be
D_{s_1\cdot \bs m,(s_1)_k\bs l,k}=x^{-m_1-1}D_{\bs m,\bs l,k}x^{m_1+1}.
\ee
The proposition follows from formula \Ref{ui}.
\end{proof}

Let  $P(\bs m,\bs l,k)$ be the monic polynomial defined by application 
of differential operators:
\bea
 P(\bs m,\bs l,k)=c(x-1)^{k+1}\tilde D_r(\bs m,\bs l,k)\tilde D_{r-1}(\bs
 m,\bs l,k)\dots \tilde D_1(\bs m,\bs l,k)\left( (x-1)^{l_1-k-1} \right) =
 \\
  c\om^{-1}(x)\
   \frac{d^{l_r-l_{r+1}}}{d x^{l_r-l_{r+1}}}\
   x^{l_r-l_{r+1}-m_r-1}\frac{d^{l_{r-1}-l_r}}{d
     x^{l_{r+1}-l_r}}\ \dots\hspace{4cm}\\ \hspace{2cm}
   x^{l_2-l_3-m_2-1}\frac{d^{l_1-l_2}}{d x^{l_{1}-l_2}}\
   \left( x^{l_1-l_2-m_1-1}(x-1)^{l_1-k-1}\right),
 \eea
where $c$ is a non-zero constant.

Following \cite{ABV}, we call $P(\bs m,\bs l,k)$
the {\it Jacobi-Pi\~neiro}
polynomial. Clearly, $P(\bs m,\bs l,k)$ is a polynomial of 
degree $l_1$ in $x$, whose
coefficients are rational functions in $\bs m, k$ with 
coefficients in $\Q$.

Sometimes we drop the dependence on $\bs m$ from our notation and simply write
$P(\bs l,k)$.
\begin{prop}\label{v=P}
Let $w\in \mc S_{r+1}$ be an element of the Weyl group, such that $w(0)=i$. 
Then
\be
v_i(\bs m, \bs l,k)=P(w\cdot \bs m, (w)_k\bs l,k).
\ee
\end{prop}
\begin{proof}
In view of Lemma \ref{v inv} it is sufficient to prove the proposition 
for $w=id$.

We have equality of functions
\bea
v_0(\bs m, l, k)= D_{\bs m, l, k}\left( 1 \right)= \hspace{7cm}\\
(x-1)^{k+1}\tilde D_r(\bs m,\bs l,k)\tilde D_{r-1}(\bs m,\bs l,k)\dots 
\tilde D_0(\bs m,\bs l,k)\ \left( \frac{1}{x-1}\right).
\eea
But
\bea
&&D_{0}(\bs m,\bs l,k)\left( \frac{1}{x-1}\right)=
\frac{d^{k-l_{1}}}
{d x^{k-l_{1}}}\left( \frac{x^{k-l_{1}}}{x-1}\right)=(k-l_1)!\ 
(x-1)^{l_{1}-k-1}.
\eea
The proposition follows.
\end{proof}

\begin{cor} Let $w\in \mc S_{r+1}$ be such that  $w(0)=0$. Then
\be
P(\bs m,\bs l,k)\ =\ P(w\cdot \bs m,(w)_k\bs l,k)\ .
\ee
\qquad $\Box$
\end{cor}

We also note that the divided Wronskian of all  elements of the basis is 
a constant, which gives an identity for Jacobi-Pi\~neiro polynomials.

\begin{cor}
We have 
\bea
 W(P(\bs m, \bs l,k),\
 P(w_1\cdot\bs m, (w_1)_k \bs l,k)\, x^{m_1+1},
 \dots,
 P(w_r\cdot \bs m, (w_r)_k\bs l,k)\, x^{\sum_{i=1}^r m_i+r})=
 \\
 {\rm const}\ x^{\sum_{i=1}^r (r+1-i)m_i}\ (x-1)^{rk}\ ,
 \eea
where $w_i = s_i \dots s_2 s_1$.
\end{cor}

\subsection{Properties of Jacobi-Pi\~neiro polynomials}
Jacobi-Pi\~neiro polynomials are examples of multiple orthogonal polynomials,
see \cite{P}, 
\cite{IN}, \cite{ABV} and references therein. In this section we prove some
properties of Jacobi-Pi\~neiro polynomials, which we need for applications.

We fix a sequence $\bs m$ of arbitrary complex numbers.

Define constants $A_i(\bs l,k)$ by
\bea
&&A_0(\bs l,k)=\prod_{s=0}^r\frac{k+\sum_{i=1}^sm_s-l_1+s+1}{k+\sum_{i=1}^sm_s-l_1+s+1+l_{s+1}-l_s}, \\
&&A_i(\bs l,k)=e_0(\bs l,k)-e_i (\bs l,k)=l_1-\sum_{s=1}^{i}m_s-k+l_{i}-l_{i+1}-i,
\eea
where $i=1,\dots,r$. 

We have
\be
A_0(\bs l,k)=(k-l_1+1)\ \prod_{s=1}^r\frac{e_s(\bs l,k)-
e_0(\bs l,k)+l_s-l_{s+1}+1}{e_s(\bs l,k)-e_0(\bs l,k)+1}.
\ee

These numbers satisfy the following property.
\begin{lem}\label{c con} We have 
\bea
A_i(\bs l+\bs 1_j,k+1)\ A_j(\bs l,k) =A_j(\bs l+\bs 1_i,k+1)\ A_i(\bs l,k).
\eea
$\qquad \Box$
\end{lem}

The recursion for spaces of polynomials 
translates to the recursion for the Jacobi-Pi\~neiro polynomials:
\begin{prop}\label{recur prop}
We have 
\be
D_i(\bs l,k)\ P(\bs l,k)=A_i(\bs l,k)\ P(\bs l+\bs 1_i,k+1).
\ee
\end{prop}
\begin{proof}
  First we suppose that 
  parameters $\bs m,\bs l,k$ and $\bs m,\bs l+\bs 1_i,k+1$ are
  consistent.
  
  By Theorem \ref {V thm} and Proposition \ref{v=P} we obtain
  that $D_i(\bs l,k) P(\bs l,k)$  is a constant multiple of 
  $P(\bs l+\bs 1_i,k+1)$. Call these constants $\tilde A_i(\bs l,k)$.
  
  Recall that $P(\bs l,k)$ is monic.  Comparing the
  highest degree terms we obtain $\tilde A_i(\bs l,k)=A_i(\bs l,k)$
  for all cases
  except for $i=0$ when the
  highest degree terms cancel and we need to consider terms of the
  next degree.

  By Lemma \ref{holon} the numbers 
  $\tilde A_i(\bs l,k)$ satisfy the relation 
  of Lemma \ref{c con}. We also have the initial conditions
\bea
\tilde A_0(\bs 0,k)=k+1.
\eea

The numbers $\tilde A_0(\bs l,k)$ are uniquely
determined by $\tilde A_i(\bs l,k)$,
$(i=1,\dots,r)$ and these properties. But numbers 
$A_0(\bs l,k)$
satisfy the same recursion. Therefore 
$\tilde A_0(\bs l,k)=\tilde A_0(\bs l,k)$.

The case of general parameters $\bs m,\bs l,k$ follows by analytic
continuation with respect to $\bs m, k$.
\end{proof}

Lemmas \ref{V 2dim} and Proposition \ref{recur prop} 
immediately imply the following 3-term 
relation for Jacobi-Pi\~neiro polynomials.
 
Let $\bs l$ be a partition. 
\begin{cor}\label{V 3 term}
We have
\be
A_i A_{js}\ P(\bs l+\bs 1_i,k+1) + A_j A_{is}\ P(\bs l+\bs 1_j,k+1)+
A_sA_{ij}\ P(\bs l+\bs 1_s,k+1)= 0,
\ee
where $A_t=A_t(\bs l,k)$, $A_{tp}=A_{tp}(\bs l,k)$. $\qquad \Box$
\end{cor}

\section{Spaces $U(\bs m,\bs l,k)$}\label{U sect}
In this section we compute explicitly a basis in the space $U(\bs m,\bs l,k)$.
\subsection{Basis in $U(\bs m,\bs l,k)$}
Let $\bs m, \bs l, k$ be consistent parameters. Recall that we use the convention $l_0=k$.

The following lemma follows from Theorem \ref{U thm}.
\begin{lem}\label{basis lem U} (Cf. \cite{Sc}.)
The space  $U(\bs m,\bs l,k)$ has a basis of the form
\be 
\{u_0(\bs m, \bs l,k),\ u_1(\bs m, \bs l,k)x^{m_r+1},\dots, 
u_r(\bs m, \bs l,k)x^{\sum_{i=1}^rm_i+r}\},
\ee
where $u_i(\bs m, \bs l,k)$ is a monic polynomial of degree
\be
\deg u_i(\bs m, \bs l,k)=l_{r-i}-l_{r-i+1}
\ee
whose coefficients are rational functions of $\bs m,k$ with coefficients in 
$\Q$.
$\qquad \Box$
\end{lem}

The polynomials $u_i$ are permuted by the Weyl group.
\begin{lem}\label{u inv}
Let $w\in \mc S_{r+1}$ be an element of the Weyl group, then
\be
u_{r-i}(w\cdot \bs m, (w)_k \bs l,k)=u_{r-w(i)}(\bs m,\bs l,k).
\ee
\end{lem}
\begin{proof}
We have the formula 
\bean\label{uP}
u_{r-i}(\bs m,\bs l,k)\ =\ c_i\ W^\dagger(f_0,f_1,\dots,\widehat{f_i},\dots,f_r)
x^{-\sum_{s=r-i+1}^{r}m_s-i},
\eean
where $c_i$ are non-zero constants, $f_i=v_i(\bs m,\bs l,k)
x^{\sum_{s=1}^{i}m_s+i}$ and $v_i$ are as in Lemma 
\ref{basis lem}. Now the lemma follows from Lemma \ref{v inv}.
\end{proof}

Note that formula \Ref{uP} determines $u_i$ in terms 
of Jacobi-Pi\~neiro polynomials.

\subsection{Recursion and 3 term relation}
Fix any sequence of $r$ complex numbers $\bs m$.
Define constants $A_i^\vee(\bs l,k)$ by
\bea
&&A_i^\vee(\bs l,k)=e_0^\vee(\bs l,k)-e_{r-i}^\vee(\bs l,k)=l_r-\sum_{s=i+1}^rm_s-l_{i}+l_{i+1}-r+i,
\\
&&A_0^\vee(\bs l,k)=-\prod_{s=0}^r\frac{\sum_{i=1}^sm_s-l_r+s}{\sum_{i=1}^sm_s-l_r+s+1+l_{r-s}-l_{r-s+1}}\ ,
\eea
where $i=0,\dots,r-1$.

These numbers satisfy the following property.
\begin{lem}\label{c con U} We have 
\bea
A_i^\vee(\bs l-\bs 1_j,k-1)\ A_j^\vee(\bs l,k) =
A_j^\vee(\bs l-\bs 1_i,k-1)\ 
A_i^\vee(\bs l,k).
\eea
In addition
\be
A_0^\vee((k,\dots,k),k)=B_{k},
\ee
where $B_{k}$ is defined by \Ref{Ak}.
$\qquad \Box$
\end{lem}

The recursion for spaces of polynomials $U(\bs m,\bs l,k)$
translates to the recursion
for the first polynomial $u_0$.
\begin{prop}\label{recur U prop}
We have 
\be
D_i^\vee(\bs l,k)\ u_0(\bs l,k)=A_i^\vee(\bs l,k)\ u_0(\bs l-\bs 1_i,k-1).
\ee

\end{prop}
\begin{proof}
Similar to the proof of Proposition \ref{recur prop}.  
\end{proof}

Let $\bs l$ be a partition such that $\bs l- \bs 1_i$, $\bs l-\bs 1_j$, 
$\bs l-\bs 1_s$ are partitions.

\begin{cor}\label{U 3 term}
We have
\be
A_i^\vee A_{js}\ u_0(\bs l-\bs 1_i,k-1) + A_j^\vee A_{is}\ 
u_0(\bs l-\bs 1_j,k-1)+
A_s^\vee A_{ij}\ u_0 (\bs l-\bs 1_s,k-1)= 0,
\ee
where $A_t^\vee=A_t^\vee(\bs l,k)$, 
$A_{tp}=A_{tp}(\bs l,k)$. $\qquad \Box$
\end{cor}
 
%Note that the relations in Corollaries \ref{V 3 term}, \ref{U 3 term}
%connect polynomials 
%with the same  $\bs m$ and $k$.
%In particular,
%the relations in Corollary \ref{V 3 term} allow us to
%obtain all polynomials $y_1(\bs l,k)$ from Jacobi polynomials
%$y_1(s,0,\dots,0)$ (cf. Section \ref{jacobi}). Similarly the relations
%in Corollary \ref{U 3 term} allow us to obtain all polynomials $y_r(\bs
%l,k)$ from polynomials $y_r(s,\dots,s)$.

\subsection{Explicit formulas}
Proposition \ref{recur U prop} allows us to compute 
the polynomials $u_i$ explicitly.

Write 
\be
u_0(\bs m,\bs l, k)=\sum_{i=0}^{l_r}(-1)^ic_i(\bs m,\bs l,k)x^{l_r-i}.
\ee
\begin{prop} \label{u exact}
We have
\be
c_i(\bs m,\bs l,k)={l_r\choose i}\prod_{j=0}^{i-1}\prod_{s=1}^r
\ \frac{\sum_{t=r-s+1}^rm_t-l_r+s+j}
{\sum_{t=r-s+1}^rm_t-l_r+s+j+1+l_{r-s}-l_{r-s+1}}.
\ee
\end{prop}
\begin{proof}
We have $c_0(\bs m,\bs l,k)=1$.
By Proposition \ref{recur U prop} we have the following relations:
\be
(l_r-i-e_{r-j}^\vee)c_i(\bs m,\bs l,k)=
A_j^\vee(\bs m,\bs l,k)c_i(\bs m,\bs l-\bs 1_j,k-1).
\ee
So, the proposition is obtained by multiplying all the constants.
\end{proof}
Note that in fact we computed all $u_i$ because of Lemma \ref{u inv}.

%----------------------------------------------------

\section{Applications to the Bethe Ansatz}\label{Bethe sec}
We consider a special case of Bethe Ansatz equation. 
It has at most one solution. If this solution exists then it is given by
zeroes of Jacobi-Pi\~neiro polynomials and of Wronskians of
Jacobi-Pi\~neiro polynomials.  
We use properties of Jacobi-Pi\~neiro polynomials to
prove that for generic values of weights, this solution
does exist and the corresponding Bethe vector is non-zero.  
However, we find that in
special cases the Bethe equation has no solutions and we obtain a
counterexample to the Bethe Ansatz Conjecture for the Gaudin model.

\subsection{Generalities of the Bethe Ansatz method}
Let $\g = \n_+\oplus\h\oplus\n_-$ be a simple Lie algebra,
\ $(.,.)$ the Killing form on $\g$,\ {}  $\al_1, \dots , \al_r \in \h^*$
simple roots. 

Let $(x_i)_{i\in I}$ be an orthonormal basis in $\g$, \
$\Omega =  \sum_{i\in I} x_i\otimes x_i\ \in \g \otimes \g$
the Casimir element. 

For a $\g$-module $V$ and $\mu \in \h^*$
denote by $V[\mu]$ the weight subspace of $V$ of weight $\mu$ and by
$\Sing\, V[\mu]$ the subspace of singular vectors of weight $\mu$,
\bea
\Sing \,V[\mu]\ =\ \{ \ v \in V\ |\ \n_+v = 0, \ hv = \langle \mu, h \rangle v \ \} \ .
\eea

Let $n$ be a positive integer and  $\bs \La = (\La_1, \dots , \La_n)$,
$\La_i \in \h^*$, a set of weights.
For a weight $\mu \in \h^*$, let $L_{\mu}$ 
be the irreducible $\g$-module with highest weight $\mu$.
Denote by $L_{\bs \La}$ the tensor product 
$L_{\La_1} \otimes \dots \otimes L_{\La_n}$.

If $X \in \End\,(L_{\La_i})$, then we denote
 by $X^{(i)} \in \End (L_{\bs \La})$ the operator
$ \cdots \otimes \id \otimes X \otimes \id \otimes \cdots$.
If $X = \sum_k X_k \otimes Y_k \in 
\End (L_{\La_i} \otimes L_{\La_j})$, then we set 
$X^{(i,j)} = \sum_k X^{(i)}_k \otimes Y^{(j)}_k\ \in \End (L_{\bs \La})$.

Let $\bs z = (z_1, \dots , z_n)$ be a point in 
$\C^n$ with distinct coordinates.
Introduce linear operators $K_1(\bs z), \dots , K_n(\bs z)$ on 
$L_{\bs \La}$ by the formula
\bea
K_i(\bs z)\ = \ \sum_{j,\ j \neq i}\ \frac{\Omega^{(i,j)}}{z_i - z_j}\ , 
\qquad i = 1, \dots , n .
\eea
The operators 
are called {\it the Gaudin Hamiltonians} of the Gaudin model associated with
$\Ve$. The  Hamiltonians commute, $ [ K_i(\bs z), K_j(\bs z) ] = 0$ 
for all $i, j$.

The problem is to diagonalize simultaneously
the  Hamiltonians, see \cite{B, BF, FFR, G, MV1, RV, ScV, V2}.

The  Hamiltonians commute with the $\g$-action on $L_{\bs \La}$, thus
is enough to
diagonalize the  Hamiltonians on the subspaces of singular vectors,\
$\Sing \,L_{\bs \La}[\mu] \subset \Ve$. 

The eigenvectors of the Gaudin Hamiltonians are constructed 
by the Bethe Ansatz method.

Fix $\La_\infty = \sum_{s=1}^n\La_s - \sum_{i=1}^r l_i\alpha_i$ where
$l_1, \dots , l_r$ are  non-negative integers. Set
\bea
\bs t\ = \ (t^{(1)}_{1}, \dots , t^{(1)}_{l_1}, t^{(2)}_{1}, \dots , t^{(2)}_{l_2},
\dots , t^{(r)}_{1}, \dots , t^{(r)}_{l_r})\ .
\eea
One defines a suitable rational function $w(\bs t,\bs z)$ with values in 
$\Ve[\La_\infty]$
as in \cite{RV}, cf. \cite{MV2, RSV}. The function is called {\it canonical}.
The canonical function  is symmetric 
with respect to the group $\bs\Sigma_{\bs l} =
\Sigma_{l_1}\times \dots \times \Sigma_{l_r}$ of permutations of coordinates $t^{(i)}_j$ 
with the same upper index. One considers the system of equations
\bean\label{Bethe eqn}
-\sum_{s=1}^n \frac{(\Lambda_s, \alpha_i)}{t_j^{(i)}-z_s}\ +\
\sum_{s,\ s\neq i}\sum_{k=1}^{l_s} \frac{(\alpha_s, \alpha_i)}{ t_j^{(i)} -t_k^{(s)}}\ +\
\sum_{s,\ s\neq j}\frac {(\alpha_i, \alpha_i)}{ t_j^{(i)} -t_s^{(i)}}
= 0, 
\eean
where $i = 1, \dots , r$ and $j = 1, \dots , l_i$.
The system is symmetric with respect to   $\bs\Sigma_{\bs l}$.
One shows that if $\bs t^0$ is a solution to \Ref{Bethe eqn}, then
$w(\bs t^0,\bs z)$ belongs to $\Sing\,\Ve[\La_\infty]$ and 
$w(\bs t^0,\bs z)$ is 
an eigenvector of
the Hamiltonians $K_1(\bs z), \dots , K_n(\bs z)$, see \cite{RV}. 

This method of finding eigenvectors is called {\it the Bethe Ansatz method}.
System \Ref{Bethe eqn} is called {\it the Bethe Ansatz equation}, 
the vector $w(\bs t^0,\bs z)$ is called {\it a Bethe vector}.

The standard form of the Bethe Ansatz Conjecture says that if
$\La_1, \dots , \La_n, \La_\infty$ are integral dominant,
 and $z_1, \dots , z_n$ are generic, then the Bethe vectors
 form a basis in $\Sing\, \Ve[\La_\infty]$. In particular, the conjecture implies
 that the number of $\bs\Sigma_{\bs l}$-orbits of solutions to \Ref{Bethe
 eqn}
 is not less than the dimension of $\Sing\, \Ve[\La_\infty]$.

\subsection{Specialized Bethe Ansatz equation}\label{bae sec}
Fix a natural number $r$ and a sequence of non-negative integers
$\bs m$. Choose a sequence of non-negative integers $\bs l$ and a
non-negative integer $k$.

Consider the following system of 
algebraic equations for variables $\bs t=(t_i^{(j)})$, 
$i=1,\dots,l_j$, $j=1,\dots,r$:
\begin{align}\label{BAE}
&\frac{m_1}{t^{(1)}_i}+\frac{k}{t^{(1)}_i-1}-
\sum_{s,\ s\neq i}\frac2{t^{(1)}_i-t^{(1)}_s}+
\sum_{s}\frac1{t^{(1)}_i-t^{(2)}_s}=0 \hspace{60pt} (i=1,\dots,l_1), \\
&\frac{m_j}{t^{(j)}_i}-\sum_{s,\ s\neq i}\frac2{t^{(j)}_i-t^{(j)}_s}+
\sum_{s}\frac1{t^{(j)}_i-t^{(j+1)}_s}+
\sum_{s}\frac1{t^{(j)}_i-t^{(j-1)}_s}=0 \quad (i=1,\dots,l_j),\notag\\
&\frac{m_r}{t^{(r)}_i}-\sum_{s,\ s\neq i}\frac2{t^{(r)}_i-t^{(r)}_s}+
\sum_{s}\frac1{t^{(r)}_i-t^{(r-1)}_s}=0\hspace{105pt}(i=1,\dots,l_r),\notag
\end{align}
where $j=2,\dots,r-1$. 

System \Ref{BAE} is the Bethe Ansatz equation \Ref{Bethe eqn}
for the Gaudin model
associated to the Lie algebra $sl_{r+1}$, points $z_1=0,z_2=1$, and weights
$\La_1,\La_2$, such that the scalar products with simple roots are given
by \Ref{weights}. The Bethe Ansatz equation with two arbitrary 
points $z_1$ and $z_2$ is related to our choice of $z_1=0$ and $z_1=1$
by a simple rescaling of variables $\bs t$.

System \Ref{BAE} is invariant with respect to the group
$\bs\Sigma_{\bs l} = \Sigma_{l_1} \times \dots \times \Sigma_{l_r}$
of permutations of variables with the same upper index. 
Thus the group acts on solutions of the system. 
This action is free and we do not distinguish 
between solutions in the same orbit.

Assume that the weight \be \La_\infty=\La_1+\La_2-\sum_s l_s\al_s \ee
is dominant. Then it is known, see \cite{MV1}, that the number of orbits of
solutions of equation \Ref{BAE} is finite and it is bounded from above
by the dimension of the space of singular vectors of weight $\La_\infty$ in the
tensor product $L_{\La_1}\otimes L_{\La_2}$.  Using Lemma \ref{decom},
we obtain the following lemma.

\begin{lem}
  If $\La_\infty$ is dominant, then system \Ref{BAE} has no solutions
  unless $\bs m, \bs l, k$ are consistent. In that case the system has
  at most one orbit of solutions.
  $\qquad \Box$
\end{lem}

%We note that solutions of the Bethe Ansatz equation 
%associated  to $sl_{r+s+1}$ with $s\in\Z_{\geq 0}$
%obviously contain solutions of the Bethe Ansatz equation
%associated to $sl_{r+1}$. Namely,
%let $\bar {\bs
%l}=(l_1,\dots,l_r,0,\dots,0)$ (with $s$ zeroes), $\bar {\bs
%m}=(m_1,\dots,m_r,n_1,\dots,n_s)$ where $n_{i}$ are any non-negative
%integers.
%\begin{lem}
%  A solution $\bs t$ of \Ref{BAE} associated to $\bs l,\bs m,k$ exists
%  if and only if there exists a solution $\bar {\bs t}$ of \Ref{BAE}
%  associated to $\bar {\bs l},\bar {\bs m},k$. Moreover if a solution
%  $\bs t$ exists and $\La_\infty$ is dominant, then $\bs t$ and
%  $\bar{\bs t}$ coincide up to permutation of variables with the same
%  upper index.
%\end{lem}
For $\bs t=(t_i^{j})$, we define a tuple of polynomials 
$\bs y=(y_1,\dots,y_r)$, 
\be
y_i=\prod_{s=1}^{l_i}(x-t_s^{(i)}).
\ee
We are interested only in zeroes of $y_i$ and we consider each 
polynomial $y_i$ up to multiplication by a non-zero complex number.
We say that $\bs y$ {\it represents} $\bs t$. Note that $\deg y_i=l_i$.

We call a tuple of polynomials $(y_1,\dots,y_r)$ {\it generic} if all
$y_i$ have
no multiple roots, no common roots with $y_{i\pm1}$ and no common
roots with $T_i$, where $T_i$ are given by \Ref{T}. If $\bs y$ represents a solution $\bs t$
of the Bethe Ansatz equation then $\bs y$ is generic.

\subsection{Connection between Bethe Ansatz and spaces of polynomials}
We fix a choice of consistent parameters $\bs m, \bs l, k$. 
Let $U=U(\bs m, \bs l,k)$ and 
$V=V(\bs m, \bs l,k)$ 
be the corresponding spaces of polynomials, see Lemmas \ref{U}, \ref{V}.

Let $f_1,\dots,f_{r+1}$ and $g_1,\dots,g_{r+1}$ be bases of $V$ and $U$
as in Lemmas \ref{basis lem} and \ref{basis lem U}. 
Then it is shown in \cite{MV1}, cf. also
Lemma \ref{dual}, that polynomials $W^\dagger(f_1,\dots,f_i)$ and
$W^\dagger(g_1,\dots,g_{r+1-i})$ are equal up to multiplication by a
non-zero constant. Let $y_i$ be the monic constant multiple
of $W^\dagger(f_1,\dots,f_i)$. We have
\bean\label{y}
y_i\ =\ c\ W^\dagger(f_1,\dots,f_i)\qquad (i=1,\dots,r),
\eean
where $c$ is a non-zero constant depending on $i$ and $f_1,\dots,f_i$.

For example in the case of $l_1=\dots=l_r=0$, we have $y_1=\dots=y_r=1$.

\begin{lem}\label{repres} 
The Bethe Ansatz equation \Ref{BAE} for consistent $\bs m, \bs l, k$
has a solution $\bs t$ if and only if tuple 
$\bs y=(y_1,\dots,y_r)$ is generic. 
Moreover, in such a case, $\bs y$ represents $\bs t$.
\end{lem}
\begin{proof}
  If the 
the Bethe Ansatz equation \Ref{BAE}, associated with consistent parameters
$\bs m, \bs l, k$,  has a solution $\bs t$, then the solution 
$\bs t$ generates ``a
  population  of solutions'', see Section 3.4 in  \cite{MV1}.
That population of solutions has its ``fundamental space'' 
$V_{\bs t}$, which is
an $r+1$-dimensional vector space of polynomials with exactly two singular 
points at 0 and 1 and exponents given by \Ref{exp1}, see Section 5.3 in
\cite{MV1}. Thus $V_{\bs t} = V(\bs m, \bs l, k)$ by Lemma \ref{V}.
Having the fundamental space $V_{\bs t}$, define the tuple $\bs y$ 
by formula \Ref{y}.
Then $\bs y$ represents the solution $\bs t$, see Theorem 5.12 in \cite{MV1}. 
Since $\bs t$ is a solution of the Bethe Ansatz equation, the tuple 
$\bs y$ is generic in the sense of Section \ref{bae sec}. 
This proves the lemma.
\end{proof}

\subsection{Properties of polynomials $y_1,\dots,y_r$}\label{y sec}
Let $\bs m,\bs l,k$ be consistent parameters. 
We have the corresponding monic polynomials  
$y_i(\bs m,\bs
l,k)$, see \Ref{y} and Lemma \ref{repres}. 

\begin{prop}\label{no zeros with T} 
For all consistent parameters $\bs m, \bs l,k$ we have
\be
y_1(\bs m,\bs l,k)(1)\neq 0, \qquad y_i(\bs m,\bs l,k)(0)\neq 0 
\qquad (i=1,\dots,r).
\ee
\end{prop}
\begin{proof}
  The proposition holds if and only if $F_{r+1-s}(0)\cap F_s(\infty)=0$ and
  $F_r(1)\cap F_1(\infty)=0$ in $V(\bs m, \bs l,k)$. The proposition is
  obviously true if $(\bs l,k)=(\bs 0,0)$. Now the proposition follows for
  all consistent parameters $\bs m, \bs l,k$ because there exists a
  canonical map $V(\bs m,\bs 0,0)\to V(\bs m, \bs l,k)$.
\end{proof}
Note that Proposition \ref{no zeros with T} shows that polynomials 
$y_i$ and $T_i$ are relatively prime.

\begin{lem}\label{rat}
The polynomials $y_i(\bs m, \bs l,k)$, $i=1,\dots,r$, 
are polynomials of degree
  $l_i$ whose  coefficients are rational functions of 
  $\bs m, k$ with coefficients in $\Q$.
\end{lem}
\begin{proof} From Theorem \ref{V thm} it is clear that we have a 
  basis in $V(\bs m,\bs l,k)$ which is compatible with the flag at
  infinity and whose coefficients are rational functions of $\bs m,k$.
  Clearly, we still have rational coefficients after the operation of
  taking the divided Wronskian.
\end{proof}

Lemma \ref{rat} allows us to consider polynomials $y_i(\bs m, \bs l,k)$ 
for all generic complex parameters $\bs m,k$.

Next two lemmas describe relations between polynomials corresponding to 
different $r$.

Let $\bar {\bs m}=(m_2,\dots,m_r)$,  $\bar {\bs l}=(l_2,\dots,l_r)$.
\begin{lem}\label{limit lem}
We have 
\bea
&&\lim_{m_1\to \infty}y_1(\bs m,\bs l, k)=(x-1)^{l_1}, \\
&&\lim_{m_1\to \infty}y_i(\bs m,\bs l, k)=
y_{i-1}(\bar{\bs m},\bar{\bs l}, l_1),
\eea
where $i=2,\dots,r$.
\end{lem}
\begin{proof}
Recall the basis of $V(\bs m, \bs l,k)$ described in Lemma \ref{basis lem}.
From \Ref{ui} we see 
that $v_i$ have finite limits when $m_1$ tends to infinity. 
Denote these limits by $\bar v_i$. From  Proposition \ref{recur prop}
we conclude that $\bar v_0=(x-1)^{l_1}$.
Consider the space
\be
\{\bar v_1,x^{m_2+1}\bar v_2,\dots,x^{\sum_{i=2}^r m_i+i-1}\bar v_r\}.
\ee
This space has the same exponents as $V(\bar {\bs m},\bar {\bs l},l_1)$.
Indeed, it is clear about exponents at infinity. 
The exponents in the finite points can only increase in the limit. 
It does not happen because of relation \Ref{reso}.

Therefore this space coincides with $V(\bar {\bs m},\bar {\bs l},l_1)$. 
The lemma follows.
\end{proof}

Let $\tilde {\bs m}=(m_1,\dots,m_r,n)$, where $n$ is any complex number, 
$\tilde {\bs l}=(l_1,\dots,l_r,0)$.
\begin{lem}
We have 
\bea
&&y_i(\tilde {\bs m},\tilde{\bs l}, k)=y_{i}(\bs m,\bs l, l_1),\\
&&y_{r+1} (\tilde {\bs m},\tilde{\bs l}, k)=1,
\eea
where $i=1,\dots,r$.
\end{lem}
\begin{proof}
The lemma follows from  Lemma \ref{V incl}.
\end{proof}

Finally we describe properties of zeros of $y_i$ in a certain asymptotic zone.

\begin{lem}\label{zone} Fix a partition $\bs l$ and a negative real number $k$.
Let $\bs m$ be a sequence of real numbers such that 
\be
m_1<<m_2<<\dots<<m_r<<0.
\ee
Then the corresponding tuple $\bs y$ is generic. 
Moreover all $y_i$ have real distinct roots $t_i^{(j)}$ 
such that $0<t_i^{(j)}<1$ and $t_{i_1}^{(j_1)}>t_{i_2}^{(j_2)}$ if $j_1<j_2$.
\end{lem}
\begin{proof}
If all $m_i<0$ then all roots of $y_1(\bs m,\bs l,k)$ are distinct, 
real and located in $(0,1)$. This is a general fact 
about roots of multiple orthogonal polynomials, see \cite{IN}.

Let $m_1<<0$, and let all other $m_i$ to be real. 
Then the roots of $y_2$ tend to the roots of another 
multiple orthogonal polynomial by Lemma \ref{limit lem}. Therefore the 
roots of $y_2$ are also real and distinct. Moreover, they are smaller
then the roots of $y_1$, which all are close to one.

Let $m_1<<m_2<<0$, and let all other $m_i$ to be real. 
Then the roots of $y_3$ also tend to the roots of a 
multiple orthogonal polynomial by Lemma \ref{limit lem}. Therefore the 
roots of $y_3$ are also real and distinct. Moreover, they are smaller
then the roots of $y_2$, which are close to one in this zone.

The lemma follows.
\end{proof}

\subsection{Bethe vectors}
%Consider the case of  $\g=sl_{r+1}$, $n=2$, $z_1=0$, $z_2=1$. Then
%$V_{\bs \La}= V_{\La_1} \otimes V_{\La_2}$. 
%Assume that for some positive integer $k$, we have
%$(\La_2,\al_1)=k$, and $(\La_2,\al_i)=0$ for $i>0$. 
Let $l_1, \dots, l_r$, $k$ be non-negative integers such that
$k\geq l_1 \geq \dots \geq l_r$. Let $m_i$ be real numbers.

%Let $m_i\geq l_i-l_{i+1}$ for all $i$. 
%Set $\mu = \La_1 + \La_2 - \sum_{i=1}^rl_i\al_i$.
%Let $t=(t^{(i)}_j)$ as before. 

We list properties of the canonical function used for the 
construction of the Bethe vectors
$w(\bs t) = w(\bs t,0,1) \in L_{\bs \La}[\La_\infty]$, see \cite{RSV}.

\begin{enumerate}
\item[(i)]
We have $L_{\bs \La}[\La_\infty] =\oplus _{\mu_1+\mu_2=\La_\infty} 
L_{\La_1}[\mu_1]\otimes
L_{\la_2}[\mu_2]$. For every such $\mu_1, \mu_2$, there exist a rational
function $w_{\mu_1,\mu_2}(\bs t) \in L_{\La_1}[\mu_1]\otimes L_{\la_2}[\mu_2]$ such that
$w(\bs t) = \sum_{\mu_1+\mu_2=\La_\infty}w_{\mu_1,\mu_2}(\bs t)$.
\item[(ii)] 
Let $F_1, \dots, F_r \in \n_-$ be the standard generators. Let 
$U(\n_-)$ be the complex universal enveloping algebra of $\n_-$. Let
$U^\R(\n_-) \subset U(\n_-)$ be its real part (consisting of polynomials in
$F_1, \dots , F_r$ with real coefficients).
For $s=1,2$, let $v_s\in L_{\La_s}$  be a highest weight vector. 
Let $L_{\La_s}^\R $ be the real part of
$L_{\La_s}$, i.e. $L_{\La_s}^\R $ is the image of $v_s$ under the action
of $U^\R(\n_-)$. For all $\mu_1, \mu_2$, we have
$w_{\mu_1,\mu_2}(\bs t) \in L_{\La_1}^\R[\mu_1]\otimes L_{\la_2}^\R[\mu_2]$
if all coordinates $t^{(i)}_j$ of the vector $\bs t$ are real.

\item[(iii)] 
Let $\bs\Sigma_{\bs l} = \prod_i \Sigma_{l_i}$ be the product of
symmetric groups.
Let $\nu$ be a function of $\bs t$. 
We define the action of $\pi \in \bs \Sigma_{\bs l}$
on $\nu$ by permuting the $t^{(i)}_j$'s with the same upper index. 
We define
the symmetrizer operator 
${\rm Sym}\,\nu    = \sum_{\pi \in \bs\Sigma_{\bs l}} \pi \cdot \nu$.
  We have the following formula 
\bea
w_{\La_1, \La_2-\sum_{i=1}^rl_i\al_i}(\bs t)\ =\hspace{8cm}\\ \hspace{3cm}
\nu(\bs t)\ v_1\otimes
F_1^{l_1-l_2}[F_2,F_1]^{l_2-l_3}\dots
[F_r,[F_{r-1},[\dots ,[F_2,F_1]\dots ]]]^{l_r}v_2\ ,
\eea
where
\bea
\nu(\bs t)\ =\ \frac 1{(l_1-l_2)!\ \dots \ (l_{r-1}-l_r)!\ l_r!}
\times \hspace{6cm}
\\
\phantom{aaa}
{\rm Sym} 
[\ \prod_{i=1}^{l_1-l_2} \frac 1{(t^{(1)}_i -1)}
\prod_{i=1}^{l_2-l_3}\frac 1
{(t^{(1)}_{i+l_1-l_2}-1)\ (t^{(2)}_{i}-t^{(1)}_{i+l_1-l_2})}
\dots \hspace{2cm}
\\
\prod_{i=1}^{l_r}\frac 1
{(t^{(1)}_{i+l_1-l_r}-1) (t^{(2)}_{i+l_2-l_r}-t^{(1)}_{i+l_1-l_r})
(t^{(3)}_{i+l_3-l_r}-t^{(2)}_{i+l_2-l_r}) \dots
(t^{(r)}_{i}-t^{(r-1)}_{i+l_{r-1}-l_r})} 
\ ]\ .
\eea
\end{enumerate}

\begin{theorem}\label{bv} Fix a partition $\bs l$. Then
for generic complex values of $k,m_1,\dots,m_r$ we have:
\begin{itemize}
\item system \Ref{BAE} has a unique orbit of solutions,

\item the multiplicity of this solution is equal to one,

\item the corresponding Bethe vector is non-zero.
\end{itemize}
\end{theorem}
\begin{proof}
Consider the tuple $\bs y (\bs m, \bs l, k) = (y_1(\bs m, \bs l, k), \dots,
y_r(\bs m, \bs l, k))$ associated with parameters $\bs m, \bs l, k$. By
Lemma
\ref{zone} the tuple is generic for generic complex values of the
parameters
 $\bs m, k$. Hence for almost all consistent values of
 parameters $\bs m, \bs l, k$
 with fixed $\bs l$, the tuple $\bs y (\bs m, \bs l, k)$ is generic.
 The exceptions lie in a proper algebraic subset in the space of all $\bs
 m, k$.

 Hence by Lemma 6.1
 for almost all consistent $\bs m, \bs l, k$ with fixed $\bs l$ the Bethe
 Ansatz equation has exactly one orbit of solutions. Equisingularity
 theorems in \cite{V1} imply that
  the Bethe Ansatz equation has exactly one orbit for almost
 all complex values of parameters
 $\bs m, k$ with fixed $\bs l$.

 Let us prove that the solutions of that only orbit have multiplicity
 one. Indeed if $\bs m, \bs l, k$ are as in Lemma \ref{zone}, then
 using the properties of the weight function described above,
 we obtain that the corresponding Bethe vector $w(t)$
 is non-zero and belongs to the real part of $L_{\bs \La}$.
 Then the multiplicity of solutions of the only orbit is equal to one
 by Corollary 7.2 in \cite{MV2}. Hence the multiplicity of solutions
 of the only orbit is equal to one
 for almost all values of $\bs m, k$. On the other hand if the multiplicity
 is one, then the corresponding Bethe vector is non-zero, see \cite{MV2}.
%For consistent values of parameters $\bs m, \bs l, k$, 
%system \Ref{BAE} has a unique orbit of solutions, 
%therefore this system has a unique orbit of solutions for
%generic values of $\bs m, k$.

%Using Lemma \ref{zone}  and properties of the weight function 
%described above, we obtain that the Bethe vector $w(t)$ 
%corresponding to this solution is non-zero and belongs to
%the real part of $L_{\bs \La}$. 
%Then the multiplicity of the solution is equal to one
%due to Corollary 7.2 in \cite{MV2}.
\end{proof}

\subsection{A counterexample to the Bethe Ansatz Conjecture}
\label{counter sec}
According to Theorem \ref{bv}, for generic values of parameters, 
there is a unique solution of the Bethe Ansatz equation \Ref{BAE} 
and a unique non-zero Bethe vector. In particular, for generic 
values of $\bs m,k$ the Bethe Ansatz Conjecture holds.
In this section we give an example of parameters when this is not true.

We consider the case of $r=2$, $l_1=2$, $l_2=1$.  Then 
polynomials $y_1$ and $y_2$ are computed to be
\be
y_2=x-\frac{m_2(m_1+m_2+1)}{(m_2+2)(k+m_1+m_2)} 
\ee and 
\be
y_1=(x-1)^2+\frac{(k-1)(2k+m_2+2m_2)}{(m_1+k-1)(m_1+m_2+k)}
(x-1)+\frac{k(k-1)}{(m_2+m_1+k)(m_1+k-1)}\ .  
\ee

\begin{lem}\label{nongen} 
Suppose $(2m_1+m_2)^2+k(4m_1-m_2^2)=0$. Then we have
$y_1=y_2^2$. In particular the pair $(y_1,y_2)$ is not generic. $\qquad \Box$
\end{lem}

Let 
\be
r=2, \quad k=49, \quad m_1=2, \quad m_2=3,\quad  l_1=2,\quad l_2=1.
\ee
Define $sl_3$ weights $\La_1$, $\La_2$ by the conditions
$(\La_1,\al_i)=m_i$, $(\La_2,\al_i)=k\delta_{1i}$. Let
$\La_\infty=\La_1+\La_2-l_1\al_1-l_2\al_2$. 
Let $L_\La$ be the irreducible $sl_3$ module of highest weight $\La$.

\begin{prop}\label{counter}
   The multiplicity of
  $L_{\La_\infty}$ in $L_{\La_1}\otimes L_{\La_2}$ is $1$.  The
  corresponding Bethe equation \Ref{BAE} has no solutions.
\end{prop}
\begin{proof}
 The parameters $\bs m,\bs l,k$ are consistent, therefore the
 multiplicity of  $L_{\La_\infty}$ in $L_{\La_1}\otimes L_{\La_2}$
 is one by Lemma \ref{decom}.
The pair $y_1(\bs m,\bs
  l,k),y_2(\bs m,\bs l, k)$ is not generic by Lemma \ref{nongen}, therefore 
  system \Ref{BAE} has no solutions.
\end{proof}
Proposition \ref{counter} provides a counterexample for the standard
version of the Bethe Ansatz Conjecture for the Gaudin model.
% which states that the
%space of singular vectors in a tensor product of irreducible
%representations has a basis of Bethe vectors if points $z_i$ are
%generic.

%It also provides a counterexample to the statement in \cite{F} 
%that each 
%eigenvector of the Gaudin Hamiltonians defines a Bethe vector.
%In particular solutions of the Bethe Ansatz equation are 
%not parameterized by Miura opers in contrary to Theorem 2, 
%Corollary 3.3, and Proposition 4.10 of \cite{F}. 

\medskip

{\bf Remark.} Here is another counterexample to the standard version of
the Bethe Ansatz conjecture.

Consider the tensor product of two adjoint $sl_3$ representations.
Then the multiplicity of the trivial module in this tensor product is
$1$.  One can explicitly check that the corresponding 3-dimensional vector
space $V$ is 
\be 
V={\rm span}\left\{ (2x-1)^2,\ (x-1)^4,\ x^4 \ \right\}.
\ee 
The space $V$ is the unique 3-dimensional space of polynomials
spanned by polynomials of degrees $2,3$ and $4$ which has two finite 
singular points $0,1$ such that the 
exponents at both singular points are $0,2,4$.

Then the corresponding pair $(y_1,y_2)=((x-1/2)^2,(x-1/2)^2)$ is not
generic and by the same argument as in Proposition \ref{counter}, 
the corresponding Bethe Ansatz equation has no solutions.

\section{Appendix}\label{app sec}
\subsection{A scalar product}
We compute the scalar product of $y_1$ and $y_r$ with respect 
to the measure defined by integration with the Jacobi type weight.

Let 
\be
\om_{m,k}=(x-1)^{-k-1}x^{-m-1}.
\ee
Set
\be
|\bs m|=\sum_{i=1}^r m_i+r-1.
\ee
Then we have $\om_{|\bs m|,k}=w$, 
where $\om$ is given by \Ref{weight}.

Assume that $m$ and $k$ are negative real numbers. 
Define a scalar product of functions $f(x)$ and $g(x)$ by the formula
\be
(f(x),g(x))_{m,k}=\int_0^1f(x)g(x)\om_{m,k}(x)dx.
\ee
The Jacobi-Pi\~neiro polynomial $P(\bs m,\bs l, k)$ is the unique 
monic polynomial of degree $l_1$ which has the property
\be
(P,x^s)_{|\bs m|,k}=0
\ee
for all $s$ of the form
\be 
s=\sum_{j=r-i+1}^rm_j+t, \qquad i\in\{0,1,\dots,r-1\},\ 
t-i\in\{0,1,\dots,l_{r-i}-l_{r-i+1}-1\},
\ee
where $k$ is a negative number and $m_j$ are such that 
$\sum_{j=0}^i{m_j}+i-1<0$ for $i=1,\dots,r$.

Operators $D_i$ and $D_i^\vee$ are adjoint in the following sense:
\begin{lem}\label{ad}
We have
\be
(D_i(\bs m,\bs l,k)f,g)_{|\bs m|+1,k+1}=
-(f,D_i^\vee(\bs m,\bs l+\bs 1_i,k+1)g)_{|\bs m|+1,k}.
\ee
\end{lem}
\begin{proof}
The lemma is checked by integration by parts.
\end{proof}

Using Lemma \ref{ad} and Propositions \ref{recur prop} and \ref{recur U prop}, 
the following scalar product is easy to compute explicitly:
\begin{lem} We have
\bea
(y_1(\bs m,\bs l,k),y_r(\bs m,\bs l,k))_{| \bs m|+1,k}= (1,1)_{| \bs m|+1,k-l_1}\times  \prod_{i=0}^{r-1} \hspace{4cm}\\
\left(\prod_{j=0}^{l_r-1}\frac{\sum_{t=r-j}^{r}m_t+i-j}
{\sum_{t=r-j}^{r}m_t+i-j+1+l_{r-i-1}-l_{r-i}}
\prod_{s=0}^{l_{i}-l_{i+1}-1}
\frac{\sum_{t=i+1}^rm_t+r+1-i+s}{\sum_{t=1}^im_t+i-s+k-l_1}\right).
\eea
$\qquad \Box$
\end{lem}

\subsection{Affine Weyl group}
In this section we note that operators $D_i$ together 
with the Weyl group generate an extension of the affine Weyl group.

We call parameters $\bs m,\bs l,k$ {\it admissible} if the polynomials
$v_i(\bs m,\bs l,k)$, $i=0,\dots,r$, are well defined. 
For admissible parameters $\bs m,\bs l,k$, we define 
the corresponding space $V(\bs m,\bs l, k)$
as follows. The space $V(\bs m,\bs l, k)$ is the space of functions with 
the explicit linear basis 
$\{v_0(\bs m,\bs l, k),v_1(\bs m,\bs l, k)x^{m_1+1},\dots, v_r(\bs m,\bs l,k)x^{\sum_{j=1}^rm_j+r}\}$, 
where $v_i(\bs m,\bs l, k)$ are given by \Ref{ui}.

Let $V$ be a disjoint  direct sum of spaces $V(\bs m,\bs l,k)$,
\be
V=\oplus_{\bs m, \bs l, k} V(\bs m,\bs l, k)
\ee
where the sum is taken over all 
sequences $\bs m$ of $r$ complex numbers, all sequences 
$\bs l$ of non-negative integers and all non-negative integers 
$k$ such that $\bs m,\bs l,k$ are admissible and
$k\geq l_1\geq \dots\geq l_r$.

The Weyl group $\mc S_{r+1}$ acts on $V$ by acting on parameters. Namely,
we define the action of the Weyl group on $V$ by 
the following rule:
\begin{align*}
s_i:\  V(\bs m,\bs l, k)&\to  V(s_i\cdot \bs m,(s_i)_k\bs l, k),\\
 v_j(\bs m,\bs l, k)x^{\sum_{t=1}^jm_t+j}&\mapsto
v_{s_i(j)}(s_i\cdot \bs m,(s_i)_k\bs l, k)x^{\sum_{t=1}^{s_i(j)}m_t+s_i(j)},
\end{align*}
where $i=1,\dots,r$ and $j=0,\dots,r$.

Introduce linear operators $D_i$, $i=0,\dots,r$, acting on $V$. 
We set
\be
D_i:\ V(\bs m,\bs l, k)\to V(\bs m, \bs l+\bs 1_i,k+1)
\ee
so that the restriction of operator $D_i$ to $V(\bs m,\bs l, k)$  equals to
$D_i(\bs m,\bs l,k)$.

We have the following commutation relations.
\begin{lem}
We have
\bean\label{rel}
D_iD_j=D_jD_i, \qquad s_iD_j=D_{s_i(j)}s_i.
\eean
\qquad $\Box$
\end{lem}

Consider the semigroup $G$ generated by $s_i$, $i=1,\dots,r$, and $D_i$, 
$i=0,\dots,r$, subject to the Weyl group relations $s_i^2=id$, 
$s_is_{i+1}s_i=s_{i+1}s_is_{i+1}$, $s_is_j=s_js_i$ ($|i-j|>1$) and 
relations \Ref{rel}.

%Let $\hat S_{r+1}$ is the affine Weyl group. It is generated by reflections 
%$\hat s_i$, $i=0,\dots,r$, subject to relations $s_i^2=id$, 
%$s_is_{i+1}s_i=s_{i+1}s_is_{i+1}$, $s_is_j=s_js_i$ ($|i-j|>1$), 
%where all indices are considered modulo $r+1$.

Note that the element $D_0\dots D_r$ is in the center of $G$. 
The factor of $G$ by the sub-semigroup generated by this 
element is isomorphic to the affine 
Weyl group of $sl_{r+1}$,
cf. \cite{K}, Proposition 6.5.

%The factor of $G$ by this element 
%is a group which is known as affine Weyl group,
%cf. \cite{K}, Proposition 6.5.
%\begin{lem}
%5We have $G/(D_0\dots D_r) \ \simeq \ \hat S_{r+1}$. $\qquad \Box$
%\end{lem}

\end{document}